\documentclass[a4paper]{amsart}
\usepackage[latin1]{inputenc}
\usepackage[all]{xy}
\usepackage{amssymb}
\usepackage{rotating}
\usepackage{amsmath}
\usepackage{graphicx}
\usepackage{epsfig}
\usepackage{eufrak}

\newtheorem{theorem}{Theorem}[section]
\newtheorem{rem}[theorem]{Remark}
\newtheorem{prop}[theorem]{Proposition}
\newtheorem{lemma}[theorem]{Lemma}
\newtheorem{defi}{Definition}[section]

\newcommand{\bprf}{{\it Proof.~}}
\newcommand{\binf}{{\it In fact }}
\newcommand{\ra}{\rightarrow}
\newcommand{\eprf}{\hfill $\square$ \smallskip\par}
\newcommand{\erem}{\hfill $\square$}

\newcommand{\C }{ \mathbb{C}}
\newcommand{\R}{ \mathbb{R}}
\newcommand{\Z}{\mathbb{Z}}
\newcommand{\Q}{\mathbb{Q}}
\newcommand{\K}{\mathbb{K}}
\newcommand{\N}{\mathbb{N}}

\newcommand{\Dh}{\mathcal{D}}
\newcommand{\Ag}{\frak{A}}

\newcounter{HPdisc}[section]
\renewcommand{\theHPdisc}{(P
)}

\newcounter{HPcont}[section]
\renewcommand{\theHPcont}{(P
)}
\newenvironment{assumeP}{
\begin{itemize}\refstepcounter{HPcont} \item[\!\!\!\!\!\!\!\!\!\!\textbf{\theHPcont}]} {\end{itemize}\par}

\input epsf
\input cyracc.def

\makeatletter
\def\blfootnote{\xdef\@thefnmark{}\@footnotetext}

\begin{document}

\title{Elliptic K3 surfaces with abelian and dihedral groups of symplectic automorphisms}
\author{Alice Garbagnati}
\address{Alice Garbagnati, Dipartimento di Matematica, Universit\`a di Milano,
  via Saldini 50, I-20133 Milano, Italia}

\email{alice.garbagnati@unimi.it}

\begin{abstract}
We analyze K3 surfaces admitting an elliptic fibration
$\mathcal{E}$ and a finite group $G$ of symplectic automorphisms
preserving this elliptic fibration. We construct the quotient
elliptic fibration $\mathcal{E}/G$ comparing its properties to the
ones of $\mathcal{E}$.\\
We show that if $\mathcal{E}$ admits an $n$-torsion section, its
quotient by the group of automorphisms induced by this section
admits again an $n$-torsion section.\\
Considering automorphisms coming from the base of the fibration,
we can describe the Mordell--Weil lattice of a fibration described
by Kloosterman. We give the isometries between lattices described
by the author and Sarti and lattices described by Shioda and by
Griees and Lam. Moreover we show that for certain groups $H$ of
$G$, $H$ subgroups of $G$, a K3 surface which admits $H$ as group
of symplectic automorphisms actually admits $G$ as group of
symplectic automorphisms.
\end{abstract}

\maketitle

\markboth{ELLIPTIC K3 SURFACES AND SYMPLECTIC AUTOMORPHISMS}{ALICE
GARBAGNATI}

\section{Introduction}
An automorphism on a K3 surface is called symplectic if it leaves
invariant the nowhere vanishing holomorphic two form of the K3
surface. Nikulin was the first to study this kind of
automorphisms. One of the main properties of a finite group $G$ of
symplectic automorphisms on a K3 surface $X$ is that the
desingularization of the quotient $X/G$ is again a K3 surface. It
is a general fact (cf. \cite{Nikulin symplectic}, \cite{Whitcher})
that if $X$ is a K3 surface admitting a finite group $G$ as group
of symplectic automorphisms, the lattice
$\Omega_G:=(H^2(X,\Z)^G)^{\perp}$ is primitively embedded in
$NS(X)$. Similarly, if $Y$ is a K3 surface obtained as
desingularization of $X/G$, there exists a lattice, called $M_G$,
(which contains the curves arising from the desingularization of
$X/G$) which is primitively embedded in
$NS(Y)$.\\

Here we analyze examples of elliptic K3 surfaces admitting
symplectic automorphisms preserving the elliptic fibration, and we
describe the
quotient surfaces.\\
In \cite{symplectic prime} and \cite{symplectic not prime}
elliptic fibrations on K3 surfaces are considered to analyze some
properties of the family of K3 surfaces admitting finite abelian
groups of symplectic automorphisms. In particular, if a K3 surface
$X$ admits an elliptic fibration $\mathcal{E}$ with a torsion
section $t$, it admits a symplectic automorphism $\sigma_t$, which
is the translation by the torsion section, and which acts as the
identity on the base of the fibration. This automorphism preserves
the fibration, hence, on the desingularization of the quotient
surface $X/\sigma_t$, there is an elliptic fibration, induced by
the one on $X$. In {\it Section \ref{section: elliptic K3 surface
and autmorphism induced by sections}}, Proposition \ref{prop: NS
MW with abelina group}, we give our main result about such
elliptic fibrations: if $\mathcal{E}$ is the generic elliptic
fibration admitting a certain group $G=\langle t\rangle$
(generated by a torsion section $t$) as Mordell--Weil group, then
the elliptic fibrations $\mathcal{E}/\sigma_t$ and $\mathcal{E}$
have isomorphic Mordell--Weil groups and the same type of singular
fibers. This implies that the N\'eron--Severi groups of the
surface $X$ and $\widetilde{X/\sigma_t}$ are isometric and we
determine them. Using the fact that $\widetilde{X/\sigma_t}$ is
obtained as quotient of a K3 surface by a group of symplectic
automorphisms, we show that its N\'eron--Severi group is isometric
to $U\oplus M_G$. In \cite{symplectic prime} and \cite{symplectic
not prime} the N\'eron--Severi group of $X$ (and hence the one of
$\widetilde{X/\sigma_t}$ since they are isometric) was described
considering its relation with the lattice $\Omega_G$, and in fact
$NS(X)$ was obtained as overlattice of finite index of
$U(n)\oplus \Omega_G$, where $n=|G|$.\\
The relation between the elliptic fibration with a certain torsion
section and the elliptic fibration obtained as quotient of this by
the automorphism induced by section can be analyzed in a more
general setting. In {\it Section \ref{section: quotient of
elliptic curves}} we recall that if $E$ is an elliptic curve
defined over a field $\mathbb{K}$ with a rational $n$-torsion
point $T_n$, and there is a primitive $n$-th root of unity in
$\mathbb{K}$, then the elliptic curve $E/\langle T_n\rangle$ has
again a rational $n$-torsion point. In the particular case of the
elliptic fibration with base $\mathbb{P}^1$ (analyzed in {\it
Section} \ref{section: automorphisms induced by sections}) this
implies that the quotient of the elliptic fibration with the
automorphism induced by an $n$-torsion section, is again an
elliptic fibration
with an $n$-torsion section.\\

Until now, we considered automorphisms of elliptic fibrations with
base $\mathbb{P}^1$ (and in particular on K3 surfaces) which act
trivially on the base. But, of course, it is possible to construct
also automorphisms of the elliptic fibrations which are induced by
automorphisms of the base of the fibrations. This kind of
automorphisms is interesting because they are related to base
changes of elliptic fibrations. In {\it Section \ref{section:
automorphisms induced by automorphisms of the basis}} we construct
elliptic fibrations with base $\mathbb{P}^1$ admitting certain
automorphisms acting also on the basis and we compute the equation
of the elliptic surface obtained as quotient of the elliptic
fibrations by these automorphisms. In {\it Section \ref{section:
automorphisms on the basis of the fibration and elliptic K3
surfaces}} we come back to the case of the K3 surfaces and we
consider the automorphisms constructed in Section \ref{section:
automorphisms induced by automorphisms of the basis}, which are
symplectic in the case of K3 surface. The main result of this part
is that we can construct elliptic fibrations such that the
N\'eron--Severi group is $U\oplus \Omega_G$ for certain abelian
groups $G$ (this is in certain sense the analogue of result of
Proposition \ref{prop: NS MW with abelina group}, where we
constructed surfaces with N\'eron--Severi $U\oplus M_G$).
Moreover, this construction gives a very nice interpretation of
the lattice $\Omega_G$, which appears in this situation as the
Mordell--Weil lattice of the fibration. As a consequence, we
proved that certain lattices found by Shioda in \cite{Shioda F5}
and \cite{Shioda F5n in generale} are isometric to the lattice
$\Omega_{\Z/n\Z}$ for $n=5,6$ computed in \cite{symplectic prime}
and
\cite{symplectic not prime}.\\
Moreover, considering automorphisms which come from the base of
the fibration, we can construct some dihedral groups of symplectic
automorphisms on K3 surface. Let $\Dh_n$ be the dihedral group of
order $2n$. By \cite{Xiao}, if $\Dh_n$ appears as group of
symplectic automorphisms of a K3 surface, then $n=3,4,5,6$. In
{\it Section \ref{section: a K3 surface with the dihedral group on
four elements as group of symplectic automorphisms}} a family of
elliptic K3 surfaces admitting the group $\Dh_4$ as group of
symplectic automorphisms is described. The automorphisms
generating $\Dh_4$ came from automorphisms of the base of the
fibration. The family considered here was described in
\cite{Kloosterman} to give an example of an elliptic K3 surface
with rank of Mordell--Weil group equal to 15. However in
\cite{Kloosterman} the Mordell--Weil lattice and the
N\'eron--Severi group of this fibration are not described. Here we
prove that the Mordell--Weil lattice has to be equal to
$\Omega_{\Dh_4}$ and, considering a special member of the family,
we compute explicitly this lattice and we prove that it is
isometric to a lattice described in \cite{EE8- lattices and
dihedral}. Moreover, (as in case of the abelian group $G$) we
prove that the N\'eron--Severi group has to
be isometric to $U\oplus \Omega_{\Dh_4}$. \\

In {\it Section \ref{section: dihedral groups of automorphisms on
elliptic fibration with torsion}} we consider examples of elliptic
K3 surfaces admitting the dihedral groups $\Dh_n$ as groups of
symplectic automorphisms. We consider both automorphisms which
come from automorphisms of the base, and which are induced by
torsion section, so in a certain sense we combine the
constructions considered before. The examples considered in this
section are used to prove an interesting phenomenon, which was
already described in \cite{dihedral 5} for the group of symplectic
automorphisms $\Dh_5$: let $H$ be a subgroup of $G$. If $G$
appears as group of symplectic automorphisms on a K3 surface, of
course so does $H$. For certain groups $H$ and $G$, if one
constructs the family of the K3 surfaces admitting $H$ as group of
symplectic automorphisms, one finds that all the members of this
family admit the group $G$ as group of symplectic automorphisms.
We prove that this happens for the pairs $(H,G)$ $(\Z/n\Z,\Dh_n)$,
$n=5,6$, $(\Z/2\Z\times \Z/4\Z,\Z/2\Z\times \Dh_4)$,
$(\Z/3\Z\times\Z/3\Z,\Ag_{3,3})$. In particular this implies that
for each of these pairs $\Omega_H\simeq \Omega_G$. In {\it Section
}\ref{section: final remarks on finite groups of symplectic
automorphisms} we reconsider the problem analyzed in Section
\ref{section: dihedral groups of automorphisms on elliptic
fibration with torsion} in a more general setting. We consider
pairs of groups $(H,G)$ such that $H$ is a subgroup of $G$ and
such that the K3 surfaces admitting $H$ as group of symplectic
automorphisms admits in fact the group $G$ as groups of symplectic
automorphisms and we show that this property is equivalent to a
property related only to the lattices $\Omega_H$ and $\Omega_G$
(Proposition \ref{prop: pairs of group (H,G)}). We prove that
there are pairs $(H,G)$ with this property which are not among the
ones described in Section \ref{section: dihedral groups of
automorphisms on elliptic fibration with torsion}, but the ones
described in Section \ref{section: dihedral groups of
automorphisms on elliptic fibration with torsion} are the only
ones such that rk$\Omega_H$=16.\\

{\it I would like to thank Bert van Geemen for his invaluable
help, Jaap Top for useful discussions and Tetsuji Shioda for the
lectures he held in Milan, which aided me greatly in working on
this paper.}

\section{Background material}
Here we collect the main results on symplectic automorphisms on K3
surfaces and on elliptic fibrations on $\mathbb{P}^1$ which will
be useful in the following.
\subsection{Symplectic automorphisms on K3
surfaces.}\label{subsection: symplectic automorphisms on K3
surfaces}

\begin{defi} An automorphism $\sigma$ on a K3 surface $X$ is symplectic
if $\sigma^*$ acts as the identity on $H^{2,0}(X)$, that is
$\sigma^*(\omega_X)=\omega_X$, where $\omega_X$ is a nowhere
vanishing holomorphic two form on $X$.\\ Equivalently $\sigma$ is
symplectic if the isometry induced by $\sigma^*$ on the
transcendental lattice is the identity.
\end{defi}

Let $\sigma$ be an automorphism of finite order on a K3 surface.
The desingularization of $X/\sigma$ is a K3 surface if and only if
$\sigma$ is symplectic.\\
The main results about the finite abelian groups of symplectic
automorphisms on K3 surfaces were obtained by Nikulin in
\cite{Nikulin symplectic}. In particular he classifies these
groups and proves that the isometries induced on the second
cohomology group of the K3 surfaces, called $\Lambda_{K3}$, by a
finite abelian group of symplectic automorphisms on a K3 surface
are essentially unique (\cite[Definition 4.6, Theorem 4.7]{Nikulin
symplectic}). As a consequence of this theorem one obtains that if
a K3 surface $X$ admits a certain group $G$ as group of symplectic
automorphisms, then the lattice
$\Omega_G:=(\Lambda_{K3}^G)^{\perp}$, which is primitively
embedded in $NS(X)\subset H^2(X,\Z)$ (cf. \cite{Nikulin
symplectic}, \cite{symplectic prime}, \cite{symplectic not
prime}), does not depend on $X$, but only on $G$. The lattices
$\Omega_G$ are computed in \cite{symplectic prime},
\cite{symplectic
not prime} for each possible finite abelian group $G$ of symplectic automorphisms on a K3 surface.\\
Hence, the main result about the finite abelian groups of
symplectic automorphisms is:
\begin{prop}{\rm (\cite{Nikulin symplectic}, \cite{symplectic prime}, \cite{symplectic not prime})}\label{prop: X has G as
symplectic group iff omegaG in NS(X)} Let $G$ be a finite abelian
group. A K3 surface $X$ admits $G$ as group of symplectic
automorphisms if and only if $\Omega_G$ is primitively embedded in
$NS(X)$.\end{prop} In \cite{Nikulin symplectic} the number and the
kind of the singularities of the quotient of a K3 surface by a
finite abelian group of symplectic automorphisms are classified.
If $X$ is a K3 surface admitting a finite abelian group of
symplectic automorphisms $G$, then $X/G$ has singularities of type
$A_l$. Hence in the desingularization $\widetilde{X/G}$ we obtain
some rational curves $M_i$ which are the exceptional curves of the
blow ups which desingularize the surface.
\begin{defi}\label{defi: MG}
Let $M_G$ be the minimal primitive sublattice of
$NS(\widetilde{X/G})$ containing the curves $M_j$ arising from the
resolution of the singularities of $X/G$. \end{defi} The lattice
$M_G$ contains some linear combinations of the $M_j$ with
rational, non integer, coefficients. For each finite abelian group
$G$ which can act symplectically on a K3 surface a description of
the lattice $M_G$ is given in
\cite[Section 7]{Nikulin symplectic}.\\

In \cite{mukai} and in \cite{Xiao} the finite (but not necessary
abelian) group of symplectic automorphisms are classified. In
\cite{Xiao} the singularities of the surface $X/G$ (where $G$ is a
finite group of symplectic automorphisms of $X$) are described.
Moreover the main result of \cite{Nikulin symplectic} (i.e. the
uniqueness of the isometry induced on the second cohomology group
by the finite abelian group of symplectic automorphisms) holds
also in the non abelian case, \cite{Whitcher}. Hence Proposition
\ref{prop: X has G as symplectic group iff omegaG in NS(X)} holds
also for finite groups $G$, not necessary abelian. The lattice
$\Omega_G$ has not yet been explicitly determined for each $G$.

\subsection{Elliptic fibration on $\mathbb{P}^1$.}\label{subsection:
elliptic fibreation on P^1}

\begin{defi} Let $X$ be a surface and $B$ a curve. An elliptic fibration $\mathcal{E}:X\ra B$ on the surface $X$
is a fibration such that the generic fiber is a smooth curve of
genus one and such that a section $s:B\ra X$ is fixed. We call
this section the zero section.
\end{defi}

We recall that every elliptic fibration can be regarded as an
elliptic curve over the function field of the basis. We will
assume $B\simeq \mathbb{P}^1$. Under this assumption each elliptic
fibration admits a minimal Weierstrass equation of the form
$$y^2=x^3+A(\tau,\sigma)x+B(\tau,\sigma),\ \ \ A(\tau,\sigma),\ B(\tau,\sigma)\in \C[\tau,\sigma]_{hom},\ \ \deg A(\tau,\sigma)=4m,\ \deg B(\tau,\sigma)=6m
$$
for a certain $m\in \N_{>0}$ and where there exists no a
polynomial $C(\tau,\sigma)$ such that
$C(\tau,\sigma)^4|A(\tau,\sigma)$ and
$C(\tau,\sigma)^6|B(\tau,\sigma)$.\\
We use the notation $A(\tau)$ and $B(\tau)$ to indicate the
polynomial $A(\tau,1)$ and $B(\tau,1)$.\\
If $m=1$, then $X$ is a rational surface, if $m=2$, $X$ is a K3
surface, if $m>2$ $X$ the Kodaira dimension of $X$, $k(X)$, is 1.
More in general $h^{2,0}(X)=dim H^{2,0}(X)=m-1$ (and a basis of
$H^{2,0}(X)$ is given by $\tau^idx/y\wedge d\tau$,
$i=0,\ldots,m-2$).\\
Hence in particular if $\deg A(\tau)\leq 4$ and $\deg B(\tau)\leq
6$, then the
surface $X$ is rational.\\

There are finitely many singular fibers. The possible singular
fibers on an elliptic fibration are described by Kodaira (cf. for
example \cite{miranda elliptic pisa}). In particular the fibers of
type $I_n$, $n>2$ are made up of $n$ rational curves meeting as a
polygon with $n$ edges. We will call $C_0$ the irreducible
component of a reducible fiber which meets the zero section. The
irreducible components of a fiber of type $I_n$ are called $C_i$
where $C_i\cdot C_{i+1}=1$ and $i\in\Z/n\Z$. Under the assumption
$C_0\cdot s=1$, these conditions identify the components
completely once the component $C_1$ is chosen, so these conditions
identify the components up to the transformation
$C_i\leftrightarrow C_{-i}$ for all $i\in\Z/n\Z$. The components
of a reducible fiber of type $I_n$ have multiplicity one, so a
section can intersect a fiber of
type $I_n$ in any component.\\

The set of the sections of an elliptic fibration form a group (the
Mordell--Weil group), with the group law which is induced by the
one on the fibers. Let $r$ be the rank of the Mordell--Weil group
(recall that if there are no sections of infinite order then
$r=0$) and $\rho=\rho(X)$ be the Picard number of the surface $X$
and $Red$ be the set $Red=\{v\in \mathbb{P}^1\,|\, F_v\mbox{ is
reducible}\}$. Then
$$\rho(X)=rk NS(X)=r+2+\sum_{v\in Red}(m_v-1)$$ (cfr. \cite[Section
7]{shioda on mordell-weil lattice}) where $m_v$ is the number of
irreducible components of the fiber $F_v$.

\begin{defi} The {\bf trivial lattice} $Tr_X$\index{Trivial lattice} (or $Tr$) of
an elliptic fibration on a surface $X$ is the lattice generated by
the class of the fiber, the class of the zero section and the
classes of the irreducible components of the reducible fibers
which do not intersect the zero section.\end{defi} The lattice
$Tr$ admits the hyperbolic plane $U$ as sublattice and its rank is
$rk (Tr)=2+\sum_{v\in Red}(m_v-1)$.
\begin{theorem}{\rm \cite[Theorem 1.3]{shioda on mordell-weil
lattice}} The Mordell--Weil group of the elliptic fibration on the
surface $X$ is isomorphic to the quotient $NS(X)/Tr=:E(K)$.
\end{theorem}

In Section 8 of \cite{shioda on mordell-weil lattice} a pairing on
$E(K)$ is defined. The value of this pairing on a section $P$
depends only on the intersection between the section $P$ and the
reducible fibers and between $P$ and the zero section. Now we
recall the definition and the properties
of this pairing.\\
Let $E(K)_{tor}$ be the set of the torsion elements in the group
$E(K)$.
\begin{lemma}\label{lemma: height pairing}{\rm \cite[Lemma 8.1, Lemma 8.2]{shioda on mordell-weil lattice}} For
any $P\in E(K)$ there exists a unique element $\phi(P)$ in
$NS(X)\otimes
\Q$ such that: \\
i) $\phi(P)\equiv (P)$ mod $Tr\otimes \Q$ (where $(P)$ is the class of $P$ modulo $Tr\otimes \Q$)\\
ii) $\phi(P)\perp Tr$.\\
The map $\phi:E(K)\ra NS(X)\otimes \Q$ defined above is a group
homomorphism such that $Ker(\phi)=E(K)_{tor}.$
\end{lemma}
\begin{theorem}\label{theorem: height formula}{\rm \cite[Theorem 8.4]{shioda on mordell-weil
lattice}} For any $P,Q\in E(K)$ let $\langle
P,Q\rangle=-\phi(P)\cdot\phi(Q)$ (where $\cdot$ is induced on
$NS(X)\otimes\Q$ by the cup product). This defines a symmetric
bilinear pairing on $E(K)$, which induces the structure of a
positive definite lattice on
$E(K)/E(K)_{tor}$.\\
In particular if $P\in E(K)$, then $P$ is a torsion section if and
only if $\langle P,
P\rangle=0$.\\
For any $P,Q\in E(K)$ the pairing $\langle-,-\rangle$ is
$$
\begin{array}{lll}
\langle P, Q \rangle&=&\chi+P\cdot s+Q\cdot s-P\cdot Q-\sum_{v\in
Red}{\rm contr}_v(P,Q)\\
\langle P, P \rangle&=&2\chi+2(P\cdot s)-\sum_{v\in Red}{\rm
contr}_v(P)
\end{array}$$
where $\chi$ is the holomorphic Euler characteristic of the
surface and the rational numbers ${\rm contr}_v(P)$ and ${\rm
contr}_v(P,Q)$ are given in the table below
\begin{eqnarray}
\begin{array}{llllll}
&I_2&I_n&I_n^*&IV^*&III^*\\
{\rm contr}_v(P)&2/3&i(n-i)/n &\left\{\begin{array}{l} 1\mbox{ if
}
i=1\\1+n/4\mbox{ if } i=n-1\mbox{ or }i=n\end{array}\right.&4/3&3/2\\
{\rm contr}_v(P,Q)&1/3&i(n-j)/n &\left\{\begin{array}{l} 1/2\mbox{
if } i=1\\2+n/4\mbox{ if } i=n-1\mbox{ or
}i=n\end{array}\right.&2/3&-
\end{array}
\end{eqnarray}
where the numbering of the fibers is the one described before, $P$
and $Q$ meet the fiber in the component $C_i$ and $C_j$ and $i\leq
j$.\end{theorem} The pairing defined in the theorem is called the
\textbf{height pairing}. This pairing will be used to determine
the intersection of the torsion sections of the elliptic
fibrations with the irreducible components of the reducible
fibers.\\
The Mordell--Weil group equipped with the height pairing is called
Mordell--Weil lattice. We remark that if $rk(Tr)=2$ (i.e.\ there
are no reducible fibres and hence no torsion sections), the height
pairing on the Mordell--Weil group is the restriction of the cup
product on $NS(X)$ to $\phi(E(K))$. So, if $rk(Tr)=2$, the
Mordell--Weil lattice is the sublattice $\phi(E)$ of $NS(X)$.
\section{Quotient of elliptic curves by groups generated by torsion points}\label{section: quotient of elliptic curves} An
elliptic curve is a pair $(E,O)$, where $E$ is a cubic curve
defined over a field $\K$ and $O$ is a rational point on this
curve. We will assume $char\K\neq 2,3$. Each elliptic curve admits
an equation of the following form
$$y^2=f_3(x),\ \ \mbox{where }f_3(x)\in\K[x],\ \ deg(f_3)=3\mbox{
and }\Delta(f_3)\neq 0,$$ where $\Delta$ is the discriminant of a
polynomial $f$. The point $O=(0:1:0)$ with respect to the
projective coordinates $(x:y:z)$ is a point on the elliptic
curve.\\
We recall the following result:
\begin{prop}\label{prop: quotient curve has the same rational torsion} Assume that $\zeta_n\in \K$ is a primitive $n$-th root of
unity. Let $E_n$ be an elliptic curve, defined over the field
$\K$,
admitting an $n$-torsion rational point $Q$.\\
Then the elliptic curve $E_n/\langle Q\rangle$ has an $n$-torsion
rational point $R$. Moreover the elliptic curve $(E_n/\langle
Q\rangle)/\langle R\rangle$ is isomorphic to the curve
$E_n$.\end{prop} \bprf Let us consider the group $E_n[n]=\{S\in
E_n(\overline{\K})\mbox{
such that }nS=0\}$.\\
For each $\sigma\in Gal(\overline{\K}/\K)$ and $S\in E_n[n]$, also
$\sigma(S)\in E_n[n]$. So $Gal(\overline{\K}/\K)$ acts on
$E_n[n]\simeq \Z/n\Z\times \Z/n\Z$. Thus we get an homomorphism
$\rho_n:Gal(\overline{\K}/\K)\ra Aut(E_n[n])\simeq
GL_2(\Z/n\Z)$.\\
Let $Q$ and $\widetilde{R}$ be two generators of $E[n]$. By
hypothesis $Q$ is a rational point and so for each $\sigma\in
Gal(\overline{\K}/\K)$, $\sigma(Q)=Q$. Hence
$\rho_n(\sigma)=\left[\begin{array}{rr}1&\alpha_{\sigma}\\0&\beta_{\sigma}\end{array}\right].$\\
Let $e_n:E_n[n]\times E_n[n]\ra \mu_n=\{n\mbox{-th roots of
unity}\}$ be the Weil pairing. It induces an isomorphism $\wedge^2
E_n[n]\ra \mu_n$. Since $\sigma$ acts on $E_n[n]$,
$\wedge^2(\sigma)$ acts on $\wedge^2 E_n[n]$. The map
$\wedge^2(\sigma)$ is the multiplication by $det(\rho_n(\sigma))$.
As $\zeta_n\in\K$, $det(\rho_n(\sigma))=1$. This implies that
$\beta_\sigma=1$ for each $\sigma\in Gal(\overline{\K}/\K)$. So
$$\rho_n(\sigma)=\left[\begin{array}{cc}
1&\alpha_\sigma\\0&1\end{array}\right]\mbox{ and
 }\sigma(\widetilde{R})=\alpha_{\sigma} Q+\widetilde{R} \mbox{ for each }\sigma\in
Gal(\overline{\K}/\K).$$ The morphism $$\phi:E_n\ra E_n/\langle
Q\rangle$$ is defined over $\K$ and thus $R:=\phi(\widetilde{R})$
is a rational $n$-torsion point on $E/\langle Q\rangle$. The
composition $E_n\ra E_n/\langle Q\rangle \ra \left(E_n/\langle
Q\rangle\right)/\langle R\rangle$ is a multiplication by $n$ on
$E_n$.\eprf

For certain values of $n$ the curve $E_n$, which admits an
$n$-torsion rational point, the quotient map $\phi$ and of the
curve $E_n/\langle Q\rangle$ can be given explicitly. For $n=2$
these equations are very well known (see.\ e.\ g.\ \cite[Pag.
79]{Silverman tate rational points on elliptic curves}), for $n=3$
they can be found in \cite{Top elliptic three torsion}. As example
here we give a similar result for $n=4$.
\begin{prop}\label{prop: elliptic four torsion} The curve $E$ admits a 4-torsion rational point if and only if, up
to a change of coordinate, it admits an equation of the form
\begin{equation*}y^2=x(x^2+(e^2-2f)x+f^2).\end{equation*}
Let us assume that $E$ admits such an equation, then $Q=(f,ef)$ is
a 4-torsion rational point. Let $G=\langle Q\rangle$.\\
The curve $E/G$ has equation
\begin{equation*}\label{formula: quotient 4 torsion}y^2=x(x^2-2(e^2+4f)x+(e^2-4f)^2).\end{equation*} If $i\in\K$ with $i^2=-1$, then the equation becomes
$y^2=x(x^2+(h^2-2g)x+g^2)$\ with $g=4f-e^2$, $h=2ie.$ In
particular if $i\in \K$, $E/G$ admits a 4-torsion rational point,
which is the image of a
4-torsion point on $E(\bar{\K})$ not in $G$.\\
The quotient map is:
$$\phi(P)=\left\{\begin{array}{ll}\left(\frac{(f + x)^2y^2}{(f -
x)^2x^2},\frac{y(f+x)((f-x)^4-4e^2fx^2)}{x^2(f-x)^3}\right)&\mbox{
if }P=(x,y)\neq (0,0), (f,ef),
(f,-ef),O_E\\
O_{E/G}&\mbox{ if }P=(x,y)= (0,0), (f,ef),
(f,-ef),O_E\end{array}\right.$$
\end{prop}\bprf The elliptic curve $E$ has
in particular a point of order two, so has an equation of the form
$y^2=x(x^2+ax+b)$ (see e.g.\ \cite[pag.\ 79]{Silverman tate
rational points on elliptic curves}). To determine $a$ and $b$ we
observe that a smooth elliptic curve of the fibration has a point
$Q$ of order four if $Q+Q=R=(0,0)$ which is the rational point of
order 2. Geometrically this means that the tangent to the elliptic
curve through $Q$ must intersect the curve exactly in the point
$R$. Taking $Q=(e,ef)$ the equation is
$$y^2=x(x^2+(e^2-2f)x+f^2).$$ To obtain the equation of the
quotient curve $E/G$, we will consider the quotient by the point
of order two $R$, and next the quotient by the image of $Q$ on the
curve $E/\langle R\rangle$. The point $R$ corresponds to $(0,0)$.
By \cite[pag.\ 79]{Silverman tate rational points on elliptic
curves} the curve $E/\langle R\rangle$ is
$$y^2=x(x^2-2(e^2-2f)x+((e^2-2f)^2-4f^2))$$
and so
$$\begin{array}{l}
y^2=x(x-e^2)(x-e^2+4f).
\end{array}$$
This curve has three 2-torsion rational points: $(0,0)$,
$(e^2,0)$, $(e^2-4f,0)$. The 4-torsion point $Q$ is mapped to the
point $\bar{Q}=(e^2,0)$ on $E/\langle R\rangle$, which is clearly
a 2-torsion point on $E/\langle R\rangle$. We consider a
transformation which sends $\bar{Q}$ to $(0,0)$. Under this
transformation the equation becomes
$$y^2=(x+e^2)x(x+4f).$$ Its quotient by the point $(0,0)$ (i.e. by
the point $\bar{Q}$) is
$$y^2=x(x^2-2(e^2+4f)x+(e^2-4f)^2)$$
If $i\in \K$, then the equation becomes
$$y^2=x(x^2+(h^2-2g)x+g^2)\ \ \mbox{ with }\ \ g=4f-e^2,\ h=2ie$$
which is the equation of an elliptic curve with a rational
4-torsion point $(g,gh)$.\\
To describe explicitly the map we have to consider the map given
in \cite[Pag. 79]{Silverman tate rational points on elliptic
curves}, to compose them with the translation $(x,y)\mapsto
(x-e^2,y)$ on $E/\langle R\rangle$, and to compose again with the
map given in \cite[Pag. 79]{Silverman tate rational points on
elliptic curves}.\erem
\section{Automorphisms induced by sections}\label{section: automorphisms induced by sections}
Now we consider an elliptic fibration $\mathcal{E}:X\ra
\mathbb{P}^1$. We will indicate with $\mathcal{E}$ also the
equation of the fibration viewed as elliptic curve over
$k(\mathbb{P}^1)$. If the equation is in the standard Weierstrass
form, we can assume that there is a zero section which is the map
$\tau\mapsto (0:1:0)$. The other sections of the elliptic
fibration $\mathcal{E}$, if any, induce automorphisms on $X$.

\begin{defi}\label{defi: autom.induced by
section} If $t$ is a section, then we call $\sigma_t$ the
automorphism defined as translation by $t$. It acts as the
identity on the basis of the fibration
$\mathbb{P}^1$.\\
The group $G=\langle \sigma_t\rangle$ is a group of automorphisms
of $X$ which preserves the elliptic fibration $\mathcal{E}$. Let
us consider the desingularization, $\widetilde{X/G}$, of the
quotient $X/G$. Since the elliptic fibration $\mathcal{E}$ is
preserved by $G$, $\widetilde{X/G}$ admits an elliptic fibration
induced by $\mathcal{E}$. We call it $\mathcal{E}/G$.
\end{defi}

The torsion sections of the fibration $\mathcal{E}:X\ra
\mathbb{P}^1$ are torsion rational points of $\mathcal{E}$
considered as
an elliptic curve over the field $k(\mathbb{P}^1)$.\\
If the section $t$ is an $n$-torsion section, the automorphisms
$\sigma_t$ is of order $n$. In particular, if $G=tors
MW(\mathcal{E})$, then $G$ is a finite abelian group of
automorphisms on $X$.
\begin{rem}{\rm Let $t$ be a section of order $n$ on $X$.
If $X$ is a rational surface, then $2\leq n\leq 6$; if $X$ is a K3
surface, then $2\leq n\leq 8$. Indeed, if $t$ is a torsion
section, then the height formula computed on $t$ is zero (Theorem
\ref{theorem: height formula}). This implies that there are a
certain number and type of reducible fibers. The Picard number of
a rational (resp. K3) elliptic surface is 10 (resp. at most 20),
this gives a bound for the order of the possible torsion sections
on these surfaces.\erem}
\end{rem} The following lemmas describe explicitly the action of
the automorphism $\sigma_t$ on the semistable reducible fibers of
the fibration. The case of instable fibers is similar, but is not
useful in the following.
\begin{lemma}{\rm \cite[Lemma 2.2]{symplectic not prime}}\label{rem: action of sigmat non-trivial other
C0} Let $t$ be an $n$-torsion section and $F_v$ be a reducible
fiber of type $I_d$ of the fibration and let us suppose that the
section $t$ meets the fiber $F_v$ in the component $C_i$.\\ Let
$r$ be the minimal number such that
$(\sigma_t)_{|F_v}^r(C_j)=C_j$, $\forall
j\in\Z/d\Z$.\\
If $\gcd(d,n)=1$ then $r=1$, $i=0$ and $\sigma_t(C_j)=C_j$, $\forall j\in\Z/d\Z$.\\
If $\gcd(d,n)\neq 1$ and $i=0$, then $r=1$ and $\sigma_t(C_j)=C_j$, $\forall j\in\Z/d\Z$.\\
If $\gcd(d,n)\neq 1$ and $i\neq 0$, then $r=d/\gcd(d,i)$,
$r|\gcd(d,n)$, and $\sigma_t(C_j)=C_{j+i}$.\end{lemma}
\begin{lemma}\label{lemma: singular fibers of the quotients}  Let $\mathcal{E}$
be a fibration with a singular fiber $F_v$ of type $I_d$ and an
$n$-torsion section $t$ as in the previous lemma and let
$G=\langle \sigma_t\rangle$. Then the image of the fiber $F_v$ on
$\mathcal{E}/G$ is a fiber of type $I_{\frac{dn}{r^2}}$.
\end{lemma}
\bprf The singular fiber $F_v/G$ has $\frac{d}{r}$ component
meeting as a polygon (because $(\sigma_t)_{|F_v}$ identifies each
component with $r$ other components). Since $\sigma_t^{r}$ fixes
the points of intersection of the components in $F_v$ (and hence
the stabilizer of each point is $\Z/\frac{n}{r}\Z$), each vertex
of the polygon $F_v/G$ is a singular point of type
$A_{\frac{n}{r}-1}$. So in $\widetilde{X/G}$ there is a tree of
$\frac{n}{r}-1$ rational curves over such a point. Hence the fiber
corresponding to $F_v$ on $\mathcal{E}/G$ is a polygon with
$\frac{d}{r}\cdot \frac{n}{r}$ edges, i.e. a fiber of type
$I_{\frac{dn}{r^2}}.$\eprf

\begin{prop}\label{prop: E and E/G same torsion, fibration} Let $\mathcal{E}_G$ be an elliptic fibration on $\mathbb{P}^1$
admitting either $G=\Z/m\Z$ or $G=\Z/m\Z\times \Z/n\Z$ as torsion
part of the Mordell--Weil group. Then $G\subset
tors(MW(\mathcal{E}_G/G))$.
\end{prop}
\bprf This follows by Proposition \ref{prop: quotient curve has
the same rational torsion} and by the fact that $\zeta_n\in
k(\mathbb{P}^1)$ for each $n$.\eprf

The equations of elliptic fibrations $\mathcal{E}_G$ for
$G=\Z/n\Z$, $n=2,3,4,5,6,8$ are the ones given in \cite[Table
1]{symplectic not prime} with different degrees for the
polynomials with variable $\tau$. If $S$ is the surface admitting
the fibration $\mathcal{E}_G$ and $h^{2,0}(S)=m-1$, then the
degrees of the polynomial are obtained multiply the degrees given
in \cite{symplectic not prime} by $m/2$. If the number obtained is
not an integer, there exists no a surface $S$ such that
$h^{2,0}(S)=m-1$ and $S$ admits an elliptic fibration of with an
elliptic fibration with that group as torsion part of the
Mordell--Weil group. In particular the elliptic fibrations
$\mathcal{E}_G$ with $G\subset MW(\mathcal{E}_G)$ are
parameterized by the coefficients of the polynomials, and hence by
an irreducible variety.
\begin{prop}\label{prop: E and E/G same reducible fibers, general case} Let $\mathcal{E}_G$ be a generic elliptic fibration on $\mathbb{P}^1$
admitting $G=\Z/n\Z$, $n=2,3,4,5,6,8$ as torsion part of the
Mordell--Weil group (here generic means that the coefficients of
the polynomials in $\tau$ have been chosen generically). Then the
singular fibers of $\mathcal{E}_G$ and of $\mathcal{E}_G/G$ are of
the
same type.\\
Moreover $\mathcal{E}_G$ and $\mathcal{E}_G/G$ have isomorphic
N\'eron--Severi group.
\end{prop}
\bprf The type of the singular fibers of an elliptic fibration can
be completely determined by the zeros of the polynomials
$A(\tau)$, $B(\tau)$ and $\Delta(\tau)$ and hence can be directly
computed by the equation on $\mathcal{E}_G$.\\
We will consider as example the case $n=6$, the others are very
similar. The singular fibers of an elliptic fibration with a
6-torsion section $t$ are $mI_6+mI_3+mI_2+mI_1$ where $\deg
k(\tau)=\deg l(\tau)=m$ with the notation of \cite{symplectic not
prime}. We deduce the intersection between the singular fibers and
the section $t$ by the height formula (cf. Theorem \ref{theorem:
height formula}). Indeed the only possibility to have a 6-torsion
section (and hence a section with height pairing equal to zero) is
that $t$ meets all the reducible fibers in the component $C_1$. By
Lemma \ref{lemma: singular fibers of the quotients}, this is
enough to describe the elliptic fibration $\mathcal{E}_G/G$ on
$\widetilde{X/\sigma_t}$. Indeed the fibers of type $I_6$ (resp.
$I_3$, $I_2$, $I_1$) on $X$ correspond to fibers of type $I_1$
(resp. $I_2$, $I_3$, $I_6$). Hence the elliptic fibration
$\mathcal{E}_G/G$
has as singular fibers $mI_6+mI_3+mI_2+mI_1$.\\
The N\'eron--Severi group of an elliptic fibration is generated by
the class of the fiber, the classes of the sections and the
classes of the reducible components of the reducible fibers. The
hypothesis on the generality of the elliptic fibration
$\mathcal{E}_G$ guarantees that
$MW(\mathcal{E}_G)=MW(\mathcal{E}_G/G)=G$. So $\mathcal{E}_G$ and
$\mathcal{E}_G/G$ have the ``same" sections, and the ``same"
singular fibers and hence they have the ``same" N\'eron--Severi
group.\erem
\section{Elliptic K3 surfaces and automorphisms induced by
sections}\label{section: elliptic K3 surface and autmorphism
induced by sections}

In this section we focalize our attention on K3 surfaces.\\
The automorphisms $\sigma_t$ induced on an elliptic K3 surface $X$
by $n$-torsions section are symplectic. Indeed on $X$ we have the
nowhere vanishing holomorphic 2-form $dx/y\wedge d\tau$, where
$dx/y$ is a 1-form on the elliptic curve in the fiber. Since
$\sigma_t$ acts trivially on the base of the fibration and as a
translation on the fibers (hence fixing the 1-forms on it),
$dx/y\wedge d\tau$ is fixed by $\sigma_t$.
\begin{prop}\label{prop: NS MW with abelina group} Let $X_G$ be a K3 surface admitting the elliptic
fibration $\mathcal{E}$ with $tors(MW(\mathcal{E}))=G$. Then $G\subset tors(MW(\mathcal{E}/G))$.\\
Moreover, if ${X_G}$ is such that $\mathcal{E}$ is the generic
elliptic fibration with $tors(MW(\mathcal{E}))=G$, then
$NS({X_G})\simeq NS(\widetilde{{X_G}/G})\simeq U\oplus M_G$.
\end{prop}
\bprf The fact that $G\subset tors(MW(\mathcal{E}))$ is
Proposition \ref{prop: E and E/G same torsion, fibration}.\\
The possible $G$ such that there exists an elliptic fibration
$\mathcal{E}$ on a K3 surface with $tors MW(\mathcal{E})=G$ are
$G=\Z/n\Z$, $2\leq n\leq 8$, $G=(\Z/n\Z)^2$ with $n=2,3,4$,
$G=\Z/2\Z\times \Z/n\Z$ with $n=4,6$. For the group $\Z/n\Z$ with
$n=2,3,4,5,6,8$ the fact that in the generic case $NS(X_G)\simeq
NS(\widetilde{X_G/G})$ is consequence of Proposition \ref{prop: E
and E/G same reducible fibers, general case}. For the groups
$\Z/2\Z\times \Z/m\Z$ with $m=4,6$ the fact that in the generic
case $NS(X_G)\simeq NS(\widetilde{X_G/G})$ can be proved as in
Proposition \ref{prop: E and E/G same reducible fibers, general
case}, considering the equation given in \cite{symplectic not
prime}. The groups $(\Z/n\Z)^2$ are the full $n$-torsion group of
the elliptic curve $\mathcal{E}$ (which is defined over
$k(\mathbb{P}^1)$). The quotient $E/E[n]$ corresponds to a
multiplication by $n$ on the curve, and hence gives an isomorphic
curve. This of course implies that
$NS(X_G)\simeq NS(\widetilde{X_G/G})$.\\
It remains the case $G=\Z/7\Z$. There exists only one possible
elliptic fibration admitting a 7-torsion section. By
\cite{shimada} this elliptic fibration has to admit 3 fibers of
type $I_7$ and 3 fibers of type $I_1$ and it has no sections of
infinte order (because its trivial lattice has rank 20 which is
the maximal possible Picard number for a K3 surface). By Lemma
\ref{lemma: singular fibers of the quotients}, $\mathcal{E}/G$
admits 3 fibers of type $I_7$ (corresponding to the fibers of type
$I_1$ on $\mathcal{E}$) and 3 fibers of type $I_1$ (corresponding
to the fibers of type $I_7$ on $\mathcal{E}$). Moreover, by
Proposition \ref{prop: E and E/G same torsion, fibration}, such an
elliptic fibration admits a 7-torsion section. Hence, as in
Proposition \ref{prop: E and E/G same reducible fibers, general
case}, $NS(X_G)\simeq
NS(\widetilde{X_G/G})$.\\
Let us prove that $NS(\widetilde{X_G/G})\simeq U\oplus M_G$ in the
case $G=\Z/5\Z$. The proof in the other cases is essentially the
same. In this case the trivial lattice is $U\oplus A_4^{\oplus
4}$. The singular fibers of this fibration are $4I_5+4I_1$. The
elliptic fibration $\mathcal{E}/G$ has $4I_1+4I_5$ as singular
fibers (cf. Lemma \ref{lemma: singular fibers of the quotients}).
In $\mathcal{E}/G$ the components ($C_i^{(j)}$, $i,j=1,\ldots,4$)
of the fibers $I_5$ which do not meet the zero section are the
rational curves arising by the resolution of the four
singularities of type $A_4$ on ${X_G}/G$ (cf. proof of Lemma
\ref{lemma: singular fibers of the quotients}). They are the
curves $M_i$ of Definition \ref{defi: MG} and hence by
\cite{Nikulin symplectic} the class
$$v=\frac{1}{5}\left[\sum_{i=1}^2
(4C_1^{(i)}+3C_2^{(i)}+2C_3^{(i)}+C_4^{(i)})+\sum_{j=3}^4
(3C_1^{(j)}+6C_2^{(j)}+4C_3^{(j)}+2C_4^{(j)})\right]$$ is a class
in $NS(\widetilde{{X_G}/G})$. The irreducible components which do
not intersect the zero section of the reducible fibers of elliptic
fibration $\mathcal{E}/G$ and the class $v$ generate the lattice
$M_G$. They are orthogonal to the class of the fiber (because they
are linear combinations of components of a fiber) and to the class
of the zero section, hence $NS(\widetilde{{X_G}/G})\hookleftarrow
U\oplus M_G$ with finite index. By definition $M_G$ is primitive
in $NS(\widetilde{{X_G}/G})$ and since $U$ is a unimodular lattice
there are no overlattice of $U\oplus M_G$ such that $M_G$ is
primitive in them. So $NS(\widetilde{{X_G}/G})\simeq U\oplus
M_G$.\erem
\begin{rem}{\rm The class $v$ in the proof of the previous proposition generates a copy of
$\Z/5\Z$ in
$NS(\widetilde{{X_G}/G})/Tr_{\widetilde{{X_G}/G}}\simeq
MW(\widetilde{{X_G}/G})$ and this shows again that in
$\mathcal{E}/G$ there is a 5-torsion section (without using
Proposition \ref{prop: E and E/G same torsion, fibration}). The
class $t_1=s+2F-v$ is a class in $NS(\widetilde{{X_G}/G})$
corresponding to a 5-torsion section of the
fibration.\erem}\end{rem} In case $G=\Z/2\Z$ Proposition
\ref{prop: NS MW with abelina group} was already proved in
\cite[Proposition 4.2]{bert Nikulin involutions}.\\

The K3 surfaces $X_G$ are examples of K3 surfaces which both admit
a finite abelian group $G$ of symlectic automorphisms and are, at
the same time, the desingularization of the quotient of a K3
surface by the same group of automorphisms. We observe that this
is not the case for a generic K3 surface with a symplectic group
of automorphisms, i.e. in general a K3 which admits a group $G$ as
finite abelian group of symplectic automorphisms cannot be
obtained as quotient of another K3 surface by the same group of
automorphisms.

\begin{rem}{\rm Since $\widetilde{X_G/G}$ admits an elliptic
fibration with $G$ as torsion part of the Mordell--Weil group, it
admits $G$ as group of symplectic automorphisms, hence we can
consider the surface $Z_G$, the desingularization of the quotient
of $\widetilde{X_G/G}$ by $G$. The surface $Z_G$ is isomorphic to
$X_G$ by the last statement in Proposition \ref{prop: quotient
curve has the same rational torsion}.\erem}\end{rem}

\section{Automorphisms induced by automorphisms of the
basis.}\label{section: automorphisms induced by automorphisms of
the basis}

Until now we considered automorphisms which leave the base of the
elliptic fibration invariant. Now we consider automorphisms which
act also on the base of the fibration. Let us consider an
automorphism $\pi$ of $\mathbb{P}^1$. It induces the map
$\nu:(x,y,\tau)\mapsto (x,y,\pi(\tau))$. In general this is not an
automorphism of the fibration, (in general it does not preserve
the fibration), but under some conditions (which depend on $\pi$)
on the polynomials $A(\tau)$ and $B(\tau)$, it is.\\
A family of examples  can be constructed considering the
automorphisms of $\mathbb{P}^1$ $\pi_n:\tau\mapsto \zeta_n\tau$,
where $\zeta_n$ is a primitive $n$-root of the unit. So
$\nu_n:(x,y,\tau)\mapsto (x,y,\zeta_n\tau)$.\\

Let $X$ be a surface with elliptic fibration $\mathcal{E}:
y^2=x^3+A(\tau)x+B(\tau)$.\\ If $A(\pi_n(\tau))=A(\tau)$ and
$B(\pi_n(\tau))=B(\tau)$, then $\nu_n$ is an automorphism of the
$X$. The conditions on $A(\tau)$ and $B(\tau)$ imply that,
$A(\tau)$ and $B(\tau)$ are actually polynomials in $\tau^n$, so
$A(\tau)=\sum_{i=0}^{l} a_i\tau^{in}$ and $B(\tau)=\sum_{i=0}^{h}
b_i\tau^{in}$ and $\deg A(\tau)=nl$, $\deg B(\tau)=nh$. The
automorphism $\nu_n$ fixes the two fibers on the point $\tau=0$
and $\tau=\infty$ (i.e. over the fixed points for $\pi$ on
$\mathbb{P}^1$). The quotient surface, $Y=X/\nu$, is smooth, and
admits an elliptic fibration with equation
$y^2=x^3+(\sum_{i=0}^{l} a_it^{i})x+\sum_{i=0}^{h} b_it^{i}.$\\
Conversely given $Y$, $X$ is obtained from $Y$ by the base change
$\tau\mapsto t=\tau^n$ (cf. \cite[VI.4]{miranda elliptic pisa}).
If the polynomials $A(\tau)=\sum_{i=0}^{l} a_i\tau^{in}$ and
$B(\tau)=\sum_{i=0}^{h} b_i\tau^{in}$ are such that there does not
exist a polynomial $C(\tau)$ such that $C(\tau)^4|A(\tau)$,
$C(\tau)^6|B(\tau)$, then the equation obtained by the base change
is the Weierstrass form of the elliptic fibration $\mathcal{E}$ on
the surface $X$. In fact we have the following commutative diagram
$$\begin{array}{ccl}X&\ra &Y=X/\nu_n\\
\downarrow&&\downarrow\\
\mathbb{P}^1_{\tau}&\ra&
\mathbb{P}^1_t=\mathbb{P}^1_{\tau}/\nu\end{array}$$ We observe
that the degree of the polynomials in $t$ of the elliptic
fibration on $X/\nu_n$ are $l$ and $h$ and the degree of the
corresponding polynomials in $\tau$ on $X$ are $nl$ and $nh$. If
$X$ for example is an elliptic K3 surface and $n=2$ (i.e. $4<2l<9$
and $6<2h<13$), then $Y=X/\nu_n$ is a rational elliptic surface
(because $l\leq 4$ and $h\leq 6$). Hence it is possible that the
Kodaira dimensions of $X$ and $X/\nu_n$ are not the same. In
particular if $X$ is a K3 surface, $X/\nu_n$ is not. In fact the
automorphisms $\nu_n$ do not fix the two holomorphic form
$dx/y\wedge d\tau$ on $X$ and hence they are not symplectic.\\

We now construct maps $\sigma_n$ which act on $\mathbb{P}^1$ as
$\pi_n$, but which act also on the fibers in such a way that one
obtains symplectic automorphisms if the surface $X$ is a K3
surface.
\begin{defi}\label{defi: automorphisms acting on the base}Let $\sigma_n$ be the map $\sigma_n: (x,y,\tau)\mapsto (\zeta_n^2x,\zeta_n^3 y,
\zeta_n \tau)$;\\
$\mu_2$ be the map $\mu_2:(x,y,\tau)\mapsto(x,y,\frac{1}{\tau})$
and \\$\varsigma_2$ be the map $\varsigma_2:
(x,y,\tau)\mapsto(\frac{x}{\tau^{4}},-\frac{y}{\tau^{6}},
\frac{1}{\tau}).$\end{defi}

\begin{prop}\label{prop: automorphisms fixing two form}
Let $\mathcal{E}_n$ be an elliptic fibration with minimal equation
$y^2=x^3+A(\tau)x+B(\tau)$. The map $\sigma_n$ is an automorphism
of the elliptic fibration $\mathcal{E}_n$ if and only if
\begin{equation}\label{formula: condition on A and B to have
symplectic base change}A(\tau)=\zeta_n^{-4} A(\zeta_n\tau)\mbox{
and }B(\tau)=\zeta_n^{-6} B(\zeta_n\tau).\end{equation} If
\eqref{formula: condition on A and B to have symplectic base
change} holds, then the elliptic fibration $\mathcal{E}_n$ has
equation: \begin{equation}\label{formula: equation with sigma
n}y^2=x^3+ \left(\sum_{i=-[\frac{4}{n}]}^h a_i\tau^{4+ni}\right)x
+\sum_{i=-[\frac{6}{n}]}^k b_i\tau^{6+ni}\end{equation} and the
elliptic fibration $\mathcal{E}_n/\sigma_n$ has equation:
$$Y^2=X^3+\left(\sum_{i=-[\frac{4}{n}]}^h
a_iT^{4+i}\right)X+\sum_{i=-[\frac{6}{n}]}^k b_iT^{6+i}.$$ In
particular by the minimality of the equation of $\mathcal{E}_n$ it
follows that $n\leq 6$.
\end{prop}
\bprf Applying the map $\sigma_n$ to $\mathcal{E}_n$, one obtains
$y^2\zeta_n^6=x^3\zeta_n^6+A(\zeta_n\tau)x\zeta_n^2+B(\zeta_n\tau)$.
It is easy to check that it is again the equation of
$\mathcal{E}_n$ if and only if \eqref{formula: condition on A and
B to have symplectic base change} holds. The equation of the
quotient surface is computed by multiplying the equation of the
elliptic fibration $\mathcal{E}_n$ by $\tau^{6n-6}$ and
considering $X:=\tau^{2n-2}x$, $Y:=\tau^{3n-3}y$, $T:=\tau^n$ as
coordinates of the quotient. The condition on $n$ comes from the
fact that if $n\geq 6$, then $\tau^4|A(\tau)$ and
$\tau^6|B(\tau)$.\erem
\begin{rem}\label{rem: sigma n and base change}{\rm The elliptic fibration $\mathcal{E}_n$ is obtained from the
elliptic fibration $\mathcal{E}_n/\sigma_n$ by a base change of
order $n$, $T\mapsto T^n$.\\
If $\mathcal{E}_n$ is a rational (resp. K3) surface,
$\mathcal{E}_n/\sigma_n$ is a rational (resp. K3) surface. This
follows from the comparison between the degree of the polynomials
in $\tau$ (resp. $T$) defining $\mathcal{E}_n$ (resp.
$\mathcal{E}_n/\sigma_n$).\\
If $\mathcal{E}_n$ has Kodaira dimension 1 it is possible that
$\mathcal{E}_n/\sigma_n$ has a lower Kodaira dimension, this
happens for example if $h=2$, $k=4$ and $n=3$ (in this case
$\mathcal{E}_n/\sigma_n$ is a K3 surface).\erem}\end{rem} The
automorphisms $\mu_2$ and $\varsigma_2$ are essentially obtained
by $\nu_2$ and $\sigma_2$ with a change of coordinates on
$\mathbb{P}^1$.

\begin{prop}\label{prop: K3 admitting sigma2 and varsigma2}Let $\mathcal{E}$ be an elliptic fibration
$y^2=x^3+A(\tau)x+B(\tau)$ such that
\begin{equation}\label{formula: condition on A and B to have t in 1/t}A(\tau)=\tau^{4m}A(\frac{1}{\tau})\mbox{  and  }B(\tau)=\tau^{6m}B(\frac{1}{\tau}).\end{equation}
Then $\mu_2$ and $\varsigma_2$ are automorphisms of the surface,\\
$A(\tau)=\sum_{i=0}^{2m}a_i(t^i+t^{4m-i})=t^{2m}\sum_{i=0}^{2m}a_i(t^{2m-i}+\frac{1}{t^{2m-i}})$
and \\
$B(\tau)=\sum_{j=0}^{3m}b_j(t^j+t^{6m-j})=t^{3m}\sum_{j=0}^{2m}b_j(t^{3m-j}+\frac{1}{t^{3m-j}})$.\\
The quotient elliptic fibrations are respectively:
$$\begin{array}{ccl}\mathcal{E}/\mu_2:&
Y^2=&X^3+X\sum_{i=0}^{2m}a_i\left(T^{2m-i}-\sum_{k=1}^{2m-i}\binom{2m-i}{k}T^{2m-i-2k}\right)+\\&&
+\sum_{j=0}^{3m}b_j\left(T^{3m-j}-\sum_{k=1}^{3m-j}\binom{3m-j}{k}T^{3m-j-2k}\right)\\
\mathcal{E}/\varsigma_2:&
\widetilde{Y}^2=&\widetilde{X}^3+\widetilde{X}(T^2-4)^2\sum_{i=0}^{2m}a_i\left(T^{2m-i}-\sum_{k=1}^{2m-i}\binom{2m-i}{k}T^{2m-i-2k}\right)+\\&&
+(T^2-4)^3\sum_{j=0}^{3m}b_j\left(T^{3m-j}-\sum_{k=1}^{3m-j}\binom{3m-j}{k}T^{3m-j-2k}\right).\end{array}$$
\end{prop}
\bprf Under the conditions \eqref{formula: condition on A and B to
have t in 1/t} the equation of $\mathcal{E}$ is invariant under
$\mu_2$ and $\sigma_2$. The equation of the quotient surfaces are
computed setting $T=t+\frac{1}{t}$, $Y=y$, $X=x$,
$\widetilde{Y}=y\frac{(t+1)^3(t-1)^3}{t^6}$,
$\widetilde{X}=x\frac{(t+1)^2(t-1)^2}{t^4}$. We use the equalities
$t^n+\frac{1}{t^n}=T^n-\sum_{k=1}^n\binom{n}{k}T^{n-2k}$,
$\frac{t^2\pm t+1}{t^2}=T\pm1$, $\frac{(t\pm 1)^2}{t}=T\pm
2$.\erem

\section{Automorphisms on the basis of the fibration and elliptic K3
surfaces}\label{section: automorphisms on the basis of the
fibration and elliptic K3 surfaces} In this section we consider
the automorphism $\sigma_n$, Definition \ref{defi: automorphisms
acting on the base}, acting on a K3 surface. In particular we saw
(Remark \ref{rem: sigma n and base change}) that the K3 surface
$X$ admitting $\sigma_n$ as automorphism is obtained by the K3
surface $\widetilde{X/\sigma_n}$ with a base change of order $n$.
These kind of base changes are considered for example in
\cite{Shioda on spehre packing}, \cite{Shioda F5}, \cite{Shioda
F5n in generale}, \cite{Kuwata}.

\begin{rem}\label{rem: sigma_n and varsigma 2 are symplectic}{\rm If $\mathcal{E}$ is an elliptic fibration on a K3 surface $X$ admitting $\sigma_n$ (resp.
$\varsigma_2$) as automorphism, then $\sigma_n$ (resp.
$\varsigma_2$) is a symplectic automorphism on $X$. Indeed the
nowhere vanishing 2-holomorphic form $dx/y\wedge d\tau$ on $X$ is
fixed by $\sigma_n^*$ (resp. $\varsigma_2^*$).}\end{rem}
\begin{prop}\label{prop: symplectic base change} Let $2\leq n\leq 6$ and let $X_n$ be a K3 surface which admits an
elliptic fibration $\mathcal{E}_n$ satisfying the condition
\eqref{formula: condition on A and B to have symplectic base
change}. Then for a generic choice of $a_i$ and
$b_i$:\begin{itemize}\item[i)] $X$ admits rank$(\Omega_{\Z/n\Z})$
independent sections of infinite order;\item[ii)] the
N\'eron--Severi group of $X_n$ is $NS(X_n)\simeq U\oplus
\Omega_{\Z/n\Z}$; \item[iii)] the Mordell--Weil lattice of the
fibration is
$MW(\mathcal{E}_n)\simeq\Omega_{\Z/n\Z}$.\end{itemize} Let
$X_{2,2}$ be a K3 surface which admits an elliptic fibration
$\mathcal{E}_{2,2}$ satisfying the conditions
$A(\tau)=A(-\tau)=\tau^8 A(\frac{1}{\tau})$,
$B(\tau)=B(-\tau)=\tau^{12} B(\frac{1}{\tau})$ (i.e. it satisfies
the conditions \eqref{formula: condition on A and B to have
symplectic base change} of Proposition \ref{prop: automorphisms
fixing two form} with $n=2$ and \eqref{formula: condition on A and
B to have t in 1/t} of Proposition \ref{prop: K3 admitting sigma2
and varsigma2}).\\ Then for a generic choice of $a_i$ and
$b_i$:\begin{itemize}\item[i')]$X_{2,2}$ admits 12 independent
sections of infinite order;\item[ii')] the N\'eron--Severi group
of $X_{2,2}$ is $NS(X_{2,2})\simeq U\oplus
\Omega_{(\Z/2\Z)^2}$;\item[iii')] the Mordell--Weil lattice of the
fibration is
$MW(\mathcal{E}_{2,2})\simeq\Omega_{(\Z/2\Z)^2}$.\end{itemize}
\end{prop}
\bprf Let $p_n$ be the number of parameters in equation
\eqref{formula: condition on A and B to have symplectic base
change}, $m_n$ the moduli of the surface $X_n$ and $\rho_n$ the
Picard number of $X_n$.  We can act on the equation satisfying the
condition \eqref{formula: condition on A and B to have symplectic
base change} with the map $(x,y,\tau)\mapsto (\lambda^2 x,
\lambda^3 y, \tau)$ (and then divide by $\lambda$), and we can act
with the automorphisms of $\mathbb{P}^1$ fixing 0 and $\infty$
(which are the points fixed by $\pi_n$). So we have a two
dimensional family of automorphisms acting on the equation of the
fibrations, hence $m_n=p_n-2$. For each $n$ we have:
\begin{eqnarray}\label{formula: mn and rkomegan}\begin{array}{r|c|c|c|c|cc}
n&2&3&4&5&6\\
\hline m_n&10&6&4&2&2\\
rk(\Omega_{\Z/n\Z})&8&12&14&16&16
\end{array}\end{eqnarray}
We recall that $\rho_n\leq 20-m_n$, by the construction of the moduli spaces of the polarized K3 surfaces.\\
By Proposition \ref{prop: automorphisms fixing two form} and
Remark \ref{rem: sigma_n and varsigma 2 are symplectic}, the K3
surface $X_n$ admits $\sigma_n$ as a symplectic automorphism of
order $n$. Hence $\Omega_{\Z/n\Z}$ is primitively embedded in
$NS(X_n)$ (cf. \cite{symplectic not prime}). The automorphism
$\sigma_n$ fixes the class of the zero section and the class of
the fiber. These two classes span a copy of the lattice $U$ in
$NS(X_n)$, so $U\hookrightarrow NS(X_n)^{\sigma_n}$. Since
$\Omega_{\Z/n\Z}=(NS(X_n)^{\sigma_n})^{\perp}$, $U\oplus
\Omega_{\Z/n\Z}\hookrightarrow NS(X_n)$, in particular $\rho_X\geq
2+rk\Omega_{\Z/n\Z}$. The rank of $U\oplus \Omega_{\Z/n\Z}$ is
equal to $20-m_n$ for each $n=2,\ldots, 6$ (cf. table
\eqref{formula: mn and rkomegan}), hence $\rho_X=20-m_n$ and the
inclusion $U\oplus \Omega_{\Z/n\Z}\hookrightarrow NS(X_n)$ has a
finite index. Since $\Omega_{\Z/n\Z}$ is primitively embedded in
$NS(X_n)$ and $U$ is a unimodular lattice, we have $NS(X_n)\simeq
U\oplus \Omega_{\Z/n\Z}$. For a generic choice of the parameters
$a_i$ and $b_i$ in the equation of $X_n$, the elliptic fibration
has no reducible fibers (if $n\neq 5$ the singular fibers are
$24I_1$, if $n=5$ the singular fibers are $2II+20I_1$), hence
$Tr_{\mathcal{E}_n}\simeq U$ and there are no torsion sections. In
particular the Mordell--Weil lattice of $\mathcal{E}_n$ is
isometric to the sublattice $\phi(E(K))\subset NS(X_n)$ (where
$\phi$ is defined in Lemma \ref{lemma: height pairing}) and
$\phi(E(K))\simeq U\oplus \Omega_{\Z/n\Z}/U\simeq
\Omega_{\Z/n\Z}$. More explicitly, if $F$ is the class of the
fiber and $s$ the class of the zero section, for each $r\in E(K)$
the map $\phi:E(K)\ra NS(X)\otimes \Q$ sends $r$ to $r+(2-k)F-s\in
\Omega_{\Z/n\Z}\subset NS(X_n)$, where $r\cdot s=k$. The height
pairing computed on $\phi(E(K))$ coincides with the intersection
form on $NS(X_n)$ restricted to $\Omega_{\Z/n\Z}$. Hence the
Mordell--Weil lattice is isometric to $\Omega_{\Z/n\Z}$. In
particular this implies that there are rank$\Omega_{\Z/n\Z}$
independent sections of infinite order. The proof in case
$X_{2,2}$ is exactly the same.\eprf

Let $2\leq n\leq 6$ and let $Y_n$ be the desingularization of the
quotient $X_n/\sigma_n$. Then $X_n$ is obtained by a base change
of degree $n$ of $Y_n$. In the following table we give the
singular fibers of $Y_n$ in the points $\tau=0$ and $\tau=\infty$
and we describe the fibers on $\tau=0$, $\tau=\infty$ of $X_n$
(with $(X_n)_{\overline{\tau}}$ we denote the fiber of the
fibration on $X_n$ over the point $\overline{\tau}$.)
$$\begin{array}{r|c|c|c|c|cc}
n&2&3&4&5&6\\
\hline (Y_n)_0, (Y_n)_\infty&I_0^*&IV^*&III^*&II^*&II^*\\
(X_n)_0, (X_n)_\infty&y^2=x^3+ax+b&y^2=x^3+1&y^2=x^3+x&II&II\\
\end{array}$$
Since $\tau=0$ and $\tau=\infty$ are fixed points of the map
$\tau\mapsto \zeta_n\tau$, $\sigma_n$ acts as an automorphism of
the fibers $(X_n)_0$ and $(X_n)_\infty$. Indeed for $n=3,4$
$(X_n)_0$ and $(X_n)_\infty$ are exactly the elliptic curves
admitting an extra automorphism. The group law of the fiber of
type $II$ corresponds to the group $\C$.\\

In \cite{Shioda on spehre packing}, \cite{Shioda F5}, \cite{Shioda
F5n in generale} certain K3 surfaces $F^{(n)}(\alpha, \beta)$,
$n=1,\ldots, 6$ are analyzed. These are K3 surfaces admitting an
elliptic fibration with equation $y^2=x^3-3\alpha
x+(t^n+1/t^n-2\beta)$. It is clear that $F^{(n)}(\alpha, \beta)$
is obtained by a base change of order $n$ applied to
$F^{(1)}(\alpha,\beta)$. As in Remark \ref{rem: sigma n and base
change} this base change is related to the automorphism
$\sigma_n$, and in fact $F^{(1)}(\alpha,\beta)$ is the quotient of
$F^{(n)}(\alpha,\beta)$ by $\sigma_n$
(Definition \ref{defi: automorphisms acting on the base}).\\
The surfaces described in \cite{Shioda F5}, \cite{Shioda F5n in
generale}, \cite{Shioda on spehre packing}, \cite{Kuwata} are
particular members of the families described in Proposition
\ref{prop: symplectic base change}.\\

Let us consider the equation $\mathcal{E}_n$ given in
\eqref{formula: equation with sigma n} with $n=5$ (resp. $n=6$),
$h=0$, $k=1$. This is the equation of the family of K3 surface
$X_n$ in Proposition \ref{prop: symplectic base change}. Up to
projective transformations, the equation $\mathcal{E}_5$ (resp.
$\mathcal{E}_6$) becomes:
\begin{equation}\label{formula: equations e5 e6} \mathcal{E}_5:\ \
y^2=x^3+a\tau^4x+(\tau^{11}+b\tau^6+\tau^1)\ \ (\mbox{resp.
}\mathcal{E}_6:\ \
y^2=x^3+a'\tau^4x+(\tau^{12}+b'\tau^6+1))\end{equation}

\begin{rem}\label{rem: isometries Shioda lattices}{\rm The Mordell--Weil lattice of the generic
element of the family $F^{(n)}(\alpha,\beta)$ is computed in
\cite[Theorem 2.4]{Shioda F5n in generale}. For $n=5,6$ the family
described in \cite{Shioda F5n in generale} and the family
described in Proposition \ref{prop: symplectic base change} are
the same (this is clear considering the equation \eqref{formula:
equations e5 e6} of $\mathcal{E}_5$, $\mathcal{E}_6$), hence the
Mordell--Weil lattice $MW(F^{(n)}_{gen})$ of \cite[Theorem
2.4]{Shioda F5n in generale} is the lattice $\Omega_{\Z/n\Z}$, for
$n=5,6$. For $n<5$ the two dimensional family described in
\cite{Shioda F5n in generale} is a subfamily of the one described
in Proposition \ref{prop: symplectic base change}, and so the
Mordell--Weil lattices of these two families are not isometric
(moreover for $n< 5$ there are reducible fibers on
$F^{(n)}_{gen}$).\erem}\end{rem}

In the following we will analyze some surfaces admitting a
dihedral group as group of symplectic automorphisms. This is based
on the following remark.
\begin{rem}\label{rem: sigman and varsigma2 generate Dhn}{\rm The maps $\sigma_n$ and $\varsigma_2$ generate the dihedral
group on $n$ elements $\Dh_n$. So if an elliptic K3 surface admits
both of them it admits the group $\Dh_n$ as group of symplectic
automorphism.\erem}
\end{rem}
In the following we give an example of elliptic K3 surfaces
admitting both the automorphisms $\sigma_n$ and $\varsigma_2$.
Moreover we give a short proof of the fact that
$\Omega_{\Z/n\Z}\simeq \Omega_{\Dh_n}$ for $n=5,6$. The case $n=5$
is proved in \cite{dihedral 5}.  In the Section \ref{section:
dihedral groups of automorphisms on elliptic fibration with
torsion} we will consider again these isometries. Another example
of elliptic K3 surfaces admitting both the automorphisms
$\varsigma_2$ and $\sigma_n$ is considered in Section
\ref{section: a K3 surface with the dihedral group on four
elements as group of symplectic automorphisms} where $n=4$.\\

The equation of the family of elliptic K3 surfaces admitting
$\sigma_n$ $n=5,6$ as automorphism is given in \eqref{formula:
equations e5 e6}. By Proposition \ref{prop: symplectic base
change} its N\'eron--Severi group is isometric to $U\oplus
\Omega_{\Z/n\Z}$ and its Mordell--Weil lattice is isometric to
$\Omega_{\Z/n\Z}$. The elliptic fibration \eqref{formula:
equations e5 e6} admits automatically also the automorphism
$\varsigma_2$ (Proposition \ref{prop: K3 admitting sigma2 and
varsigma2}). Both automorphisms $\sigma_n$ and $\varsigma_2$ act
as the identity on the class of the fiber and on the class of the
zero section and $\sigma_n$ acts non trivially on all the other
sections. As in Proposition \ref{prop: symplectic base change},
this implies that $MW(\mathcal{E})\simeq \Omega_{\Dh_n}$, but, by
Proposition \ref{prop: symplectic base change},
$MW(\mathcal{E})\simeq \Omega_{\Z/n\Z}$ and then we conclude that
$\Omega_{\Z/n\Z}\simeq \Omega_{\Dh_n}$ for $n=5,6$.

\section{The group $\Dh_4$ as group of symplectic
automorphisms.}\label{section: a K3 surface with the dihedral
group on four elements as group of symplectic automorphisms}

Here we give an example of a family of K3 surfaces admitting the
dihedral group on four elements, $\Dh_4$, as group of symplectic
automorphisms. These automorphisms are induced on elliptic K3
surfaces by the automorphisms $\sigma_4$, $\varsigma_2$
(Definition \ref{defi: automorphisms acting on the base}) hence
they are related to a base change of the fibration (Remark
\ref{rem:
sigma n and base change}).\\
In \cite{Kloosterman} the same family of elliptic K3 surfaces is
considered in order to construct an elliptic K3 surface with 15
independent sections of infinite order. Its construction is based
on two base changes (one of order 2 and one of order 4). We prove
that the family of elliptic K3 admitting $\Dh_4$ as group of
symplectic automorphisms that we consider is the same as the one
described in \cite{Kloosterman}. As in Proposition \ref{prop:
symplectic base change}, using the fact that these surfaces admit
a certain group of symplectic automorphisms, we give a description
of the N\'eron--Severi groups and of the Mordell--Weil lattices.

\begin{rem}\label{rem: K3 with dihedral group}{\rm The elliptic K3 surfaces invariant under $\sigma_4$ and
$\varsigma_2$ have equations:
\begin{equation}\label{formula: equation elliptic with
dihedral}\mathcal{E}:\ \ \
y^2=x^3+x(a\tau^8+b\tau^4+a)+(c\tau^{10}+d\tau^{6}+c\tau^2).\end{equation}
Such an equation depends on 4 parameters, but we can use the
transformation $(x,y)\mapsto (\lambda^2x, \lambda^3y)$ to put one
of this parameters equal to 1. Hence this family has are 3 moduli.
So the Picard number of the generic K3 surface with such an
elliptic fibration is at most 17. By Remark \ref{rem: sigman and
varsigma2 generate Dhn} each member of this family admits $\Dh_4$
as group of symplectic automorphisms.\\
Let $G=\langle\sigma_4,\varsigma_2\rangle$. Then the elliptic
fibration $\mathcal{E}/G$ has equation
$$\mathcal{E}/G:\ \ Y^2=X^3+X(T^2-4)^2(aT+b)+(T^2-4)^3(cT+d),$$
where $T=(t^4+1/t^4)$, $X=x(t^2+1/t^6)$, $Y=y(t-1/t^7)$.\\
For a generic choice of $a,b,c,d$ the elliptic fibration
$\mathcal{E}/G$ has $III^*+2I_0^*+3I_1$ as singular fibers and
does not admit sections of infinite order. Indeed the trivial
lattice of this fibration has rank 17 and thus the Picard number
is at least 17. Since the Picard number of two K3 surfaces are the
same if there exists a dominant rational map between them, the
Picard number of this surface and of the one with fibration
$\mathcal{E}$ is exactly 17.\erem}\end{rem}
\begin{rem}\label{rem: kloosterman}{\rm The elliptic fibration $\mathcal{E}/G$ has the properties of the elliptic fibration
$\widetilde{\pi}:\widetilde{X}\ra \mathbb{P}^1$ considered in
\cite[Table 2]{Kloosterman}. In particular these elliptic
fibrations coincide under the condition given in \cite[Proposition
3.2]{Kloosterman} that the rank of a certain twist (in
\cite[Construction 3.1]{Kloosterman}) is 0. This implies that the
elliptic fibration $\mathcal{E}$ is the one described in
\cite[Theorem 1.1]{Kloosterman} as elliptic K3 with Mordell--Weil
rank 15. Indeed this fibration is constructed by Kloosterman
performing two base changes on $\mathcal{E}/G$ (one of order 4 and
one of order 2), which correspond to the presence of the
symplectic automorphisms $\sigma_4$ and $\varsigma_2$ on
$\mathcal{E}$.\erem}\end{rem}

In the following we construct an elliptic fibration (which is a
special member of the family $\mathcal{E}$ given in Remark
\ref{rem: K3 with dihedral group}) admitting $\Dh_4$ as group of
symplectic automorphisms and use it to determine the lattice
$\Omega_{\Dh_4}$. We will use the description of this lattice to
find the Mordell--Weil lattice and the N\'eron--Severi group of
the surface constructed in \cite{Kloosterman}.
\subsection{Kummer surface of the product of two elliptic curves}
Let us recall the construction of a Kummer surface $Km(A)$. Let
$A$ be an Abelian surface. Let $\iota_A:a\mapsto -a$. Hence
$\iota$ is an involution fixing the sixteen point of $A[2]=\{a\in
A\mbox{ such that }2a=0\}$. The quotient $A/\iota$ is a singular
surface with sixteen singularities of type $A_1$. Let us denote
with $Km(A)$ the desingularization of this surface. The surface
$Km(A)$ is a K3 surface, called Kummer surface of $A$. We have the
following commutative diagram:
$$
\begin{array}{rcl}
\iota\circlearrowright A&\stackrel{\beta}{\leftarrow}&\widetilde{A}\circlearrowleft\widetilde{\iota}\\
\pi\downarrow&&\downarrow\widetilde{\pi}\\
A/\iota&{\leftarrow}&Km(A)\simeq \widetilde{A}/\widetilde{\iota}
\end{array}
$$
$\beta$ is the blow up of $A$ in the points of $A[2]$,
$\widetilde{\iota}$ is the involution
induced by $\iota$ on $\widetilde{A}$, $\pi$ and $\widetilde{\pi}$ are the quotient maps.\\

Let $C$ and $D$ two elliptic curves with equation:
\begin{eqnarray}\label{equations elliptic curves}\begin{array}{ll}
C:\ v^2=u^3+Au+B\\
D:\ \nu^2=f(\tau)\ \ deg(f(\tau))=4,3.\end{array}
\end{eqnarray} Then (as in \cite{Kuwata}) we consider the equation
$f(\tau)y^2=x^3+Ax+B$, or, which is the same under the
transformation $y\mapsto y/f^2(\tau)$ and $x\mapsto x/f^3(\tau) $,
the equation
\begin{equation}\label{equation: kummer product elliptic curve}y^2=x^3+Af^2(\tau)x+Bf^3(\tau).\end{equation}
Let $f:\mathcal{E}\ra\mathbb{P}^1$ be the elliptic fibration with
equation \eqref{equation: kummer product elliptic curve}.\\
The fibration $\mathcal{E}$ is an isotrivial fibration with fibers
isomorphic to $C$. Indeed under the change of coordinates
\begin{eqnarray}\label{equation: changes coordinate kummer product}\left\{\begin{array}{l}y=\nu^3v\\x=\nu^2
u.\end{array}\right.\end{eqnarray} the equation \eqref{equation:
kummer product elliptic curve} becomes
$\nu^6v^2=\nu^6u^3+Af^2(\tau)\nu^2u+Bf^3(\tau)$. Since
$\nu^2=f(\tau)$, we obtain that the fiber over each fixed
$\overline{\tau}$ has equation $\nu^6(v^2=u^3+Au+B$).\\

The equation \eqref{equation: kummer product elliptic curve} is
the equation of an elliptic fibration on the Kummer surface
$Km(C\times D)$. In fact let us consider the trivial elliptic
fibration $C\times D\ra D$ defined over the Abelian surface
$C\times D$ and  the map $\varphi: D\ra \mathbb{P}^1$ defined as
$(\tau, \nu)\mapsto (\tau, -\nu)$. Then we have the following
commutative diagram:
\begin{eqnarray*}\begin{array}{cccl}\Phi:&C\times D&\stackrel{2:1}{\rightarrow}&\mathcal{E}\\
&\downarrow&&\downarrow f\\
\varphi:&D&\stackrel{2:1}{\rightarrow}&\mathbb{P}^1_{\tau}\end{array}\end{eqnarray*}
where the map $\Phi$ is induced by the map $\varphi$ and hence it
has degree 2. The explicit equation of $\Phi$ descends by
\eqref{equation: changes coordinate kummer product} and is $$\Phi:
((u,v),(\nu, \tau))\mapsto (x,y,\tau)=(\nu^2u, \nu^3v, \tau ).$$
The points $(p,q)=((u,v),(\tau, \nu))$, $(p',q')=((a,b),(\alpha,
\beta))\in C\times D$ are such that $\Phi(p,q)=\Phi(p',q')$ if and
only if either $p=p'$ and $q=q'$, or $u=a$, $v=-b$, $\tau=\alpha$,
$\nu=-\beta$. The maps $\iota_{C}:c\mapsto -c$, $\iota_D:d\mapsto
-d$ is represented in coordinates by
\begin{eqnarray*}\left\{\begin{array}{lll}(u,v)\mapsto (u,-v)&\mbox{ on the elliptic curve }&C\\
(\tau, \nu)\mapsto (\tau, -\nu)&\mbox{ on the elliptic curve
}&D\end{array}\right.\end{eqnarray*} Hence the map $\Phi$
identifies the point $(p,q)$ with $(-p,-q)$, so it coincides with
the map $\iota_{C\times D}:a\mapsto -a$ on the Abelian surface
$C\times D$. In particular, the desingularization of the quotient
$C\times D/\iota_{C\times D}$ is the surface $Km(C\times D)$ and
the singular elliptic fibration $\mathcal{E}$ induces an elliptic
fibration over $Km(C\times D)$.\\ With an abuse of notation we
will call $\mathcal{E}$ also the elliptic fibration over
$Km(C\times D)$.

\subsection{$Km(E_i\times E_i)$} From now on we assume that the two elliptic curves
$C$ and $D$ are $E_i$, i.e. they admit an automorphism of order 4
which is not a translation.
\begin{rem}\label{rem: equation Ei and automorphisms of order 4 }{\rm Since the elliptic curves $C\simeq D$ are isomorphic to
$E_i$, the equations \eqref{equations elliptic curves} can be
chosen on the following way
\begin{eqnarray*}\begin{array}{ll}
C:\ v^2=u^3+u\\
D:\ \nu^2=\tau^4-1\end{array}
\end{eqnarray*} where the
automorphisms of order 4 are $\varphi_{C}:(u,v)\mapsto (-u, iv)$
and $\varphi_{D}:(\tau, \nu)\mapsto (i\tau,\nu)$ respectively.
Hence the equation \eqref{equation: kummer product elliptic curve}
becomes \begin{equation}\label{equation: elliptic fibration
Km(Ei*Ei)}y^2=x^3+(\tau^4-1)^2x.\end{equation} The zero of the
group law on the curve $C$ (resp. $D$) has to be fixed by the map
$\iota_{C}$ (resp. $\iota_{D}$) hence we can assume that it is the
point at infinity (resp. the point
$(\tau,\nu)=(1,0)$).}\erem\end{rem} We observe that this elliptic
fibration is a special member of the family described in Remark
\ref{rem: K3 with dihedral group}, indeed its equation is obtained
by equation \eqref{formula: equation elliptic with dihedral}
putting $a=1$, $b=-2$, $c=d=0$. Hence it admits the automorphisms
$\sigma_4$ and $\varsigma_2$ and so the group $\Dh_4$ is a group
of symplectic automorphisms on it. We now relate the automorphisms
$\sigma_4$ and $\varsigma_2$ with automorphisms of
order 4 on the elliptic curves $C$ and $D$.\\

The points of $C[2]$ are, in homogeneous coordinates,
$0_C=(0:1:0)$, $c_1=(0:0:1)$, $c_2=(i:0:1)$, $c_3=(-1:0:1)$. They
correspond to the image of the point $0_C=0$, $c_1=(1+i)/2$,
$c_2=1/2$, $c_3=i/2$ under the map $\eta: \C\ra\C/\Z[i]=E_i$. They
can be identified with the points $0_C=(0,0)$, $c_1=(1,1)$,
$c_2=(1,0)$, $c_3=(0,1)$ under the identification of $C[2]$ with
$(\Z/2\Z)^2$.\\
The points of $D[2]$ are the points fixed by the map $\iota_D$.
Hence they are $0_D=(1;0)$, $d_1=(i;0)$, $d_2=(-1;0)$,
$d_3=(-i;0)$. They are the image of the points $0_D=(0,0)$,
$d_1=i/2$, $d_2=(1+i)/2$, $d_3=1/2$ under $\eta$ and they can be
identified with the points $0_D=(0,0)$, $d_1=(0,1)$, $d_2=(1,0)$,
$d_3=(1,1)$ under the identification of $D[2]$ with
$(\Z/2\Z)^2$.\\
The map $\varphi_C$ fixes the points $0_C$ and $c_1$ and switches
the point $c_2$ and $c_3$. On $\C$ it corresponds to $z_1\mapsto
iz_1$.\\
The map $\varphi_D$ acts in the following way: $0_D\mapsto
d_1\mapsto d_2\mapsto d_3$. On $\C$ it corresponds to $z_2\mapsto
iz_2+1/2$.\\

We will denote with $K_{e,f,g,h}$ the rational curve in
$Km(E_i\times E_i)$ which is the exceptional curve over the points
$p_{e,f,g,h}=(e,f)\times (g,h)\in (E_i\times E_i)[2]$.\\
Let $\Delta=\{(c,d)\in Km(E_i\times E_i)| c=d\}$ and
$\Gamma=\{(c,d)\in Km(E_i\times E_i)| c=\varphi_{E_i}(d)\}$ where
$\varphi_{E_i}$ is the automorphism of order four of $E_i$. They
correspond to two classes in N\'eron--Severi group.
\begin{rem}\label{rem: basis NS(C*D) coordinate}{\rm We identify $C$ with $\R_{x_1,x_2}^2/\Z^2$
(i.e. the real coordinates of $\R^2$ will be denoted by $x_1$ and
$x_2$), and $D$ with $\R^2_{x_3,x_4}/\Z^2$. Since $dx_i\wedge
dx_j$, $i<j$ is a basis of $H^2(C\times D,\Z)$, we can express the
classes $[C]$, $[D]$, $[\Gamma]$, $[\Delta]$ as linear combination
of $dx_i\wedge dx_j$. More precisely
$$\begin{array}{llllll} [C]&=&dx_3\wedge dx_4,\ &\left[\Gamma\right]&=&dx_1\wedge dx_2+ dx_3\wedge dx_4+ dx_1\wedge
dx_3+dx_2\wedge dx_4,\\ \left[D\right]&=&dx_1\wedge dx_2,\
&\left[\Delta\right]&=&dx_1\wedge dx_2+dx_3\wedge dx_4+dx_2\wedge
dx_3-dx_1\wedge dx_4.
\end{array}$$
This computation is done as in \cite[Remark
3.11]{alinikulin}.}\erem\end{rem}

\begin{rem}\label{rem: divisible classes in Km(CD)}{\rm Let $\omega_{i,j}=\pi_*(\beta^*(dx_i\wedge dx_j))$. By Remark \ref{rem: basis NS(C*D) coordinate}, we
have
$$\begin{array}{llllll} \pi_*(\beta^*[C])&=&\omega_{3,4},&\pi_*(\beta^*[\Gamma])&=&\omega_{1,2}+ \omega_{3,4}+ \omega_{1,3}+\omega_{2,4},\\
\pi_*(\beta^*[D])&=&\omega_{1,2},&
\pi_*(\beta^*[\Delta])&=&\omega_{1,2}+\omega_{2,3}+\omega_{3,4}-\omega_{1,4}.
\end{array}$$
Hence:
$$\begin{array}{lll} \widetilde{C_0}:=\pi(\widetilde{C\times(0,0)})&=&\frac{1}{2}\left(\omega_{3,4}-\sum_{(a,b)\in (\Z/2\Z)^2}K_{a,b,0,0}\right),\\
\widetilde{D_0}:=\pi(\widetilde{(0,0)\times D})&=&\frac{1}{2}\left(\omega_{1,2}-\sum_{(a,b)\in (\Z/2\Z)^2}K_{0,0,a,b}\right),\\
\widetilde{\Gamma}:=\pi(\widetilde{\Gamma})&=&\frac{1}{2}\left(\omega_{1,2}+ \omega_{3,4}+ \omega_{1,3}+\omega_{2,4}-\sum_{(a,b)\in (\Z/2\Z)^2}K_{a,b,b,a}\right),\\
\widetilde{\Delta}:=\pi(\widetilde{\Delta})&=&\frac{1}{2}\left(\omega_{1,2}+\omega_{2,3}+\omega_{3,4}-\omega_{1,4}-\sum_{(a,b)\in
(\Z/2\Z)^2}K_{a,b,a,b}\right).
\end{array}$$
\erem}
\end{rem}

We use the following convention: if $W$ is a vector subspace of
the affine space $(C\times D)[2]=(\Z/2\Z)^4$, then $\bar{K}_W$ is
the class $\frac{1}{2}\sum_{a\in W}K_{a}$ and
$\hat{K}=\frac{1}{2}\sum_{a\in (\Z/2\Z)^4}K_{a}=\bar{K}_{A[2]}$.
Moreover, $W_i=\{a=(a_1,a_2,a_3,a_4)\in A[2]\mbox{ such that
}a_i=0\}$, $i=1,2,3,4$.

\begin{prop}\label{prop: action of symplectic dihedral} Let $Km(C\times D)$ be the K3 surface admitting the elliptic fibration
$$y^2=x^3+(\tau^4-1)^2x.$$
A $\Z$-basis for $NS(Km(C\times D))$ is given by the classes:
$$\begin{array}{cccccccccc}
K_{0000},&K_{0001},&K_{0010},&K_{0100},&K_{1000},&K_{0011},&K_{0101},&K_{1001},&K_{0110},&K_{1010},\\
K_{1100},&\bar{K}_{W_4},&\bar{K}_{W_3},&\bar{K}_{W_2},&\bar{K}_{W_1},&\hat{K},&
\widetilde{C_0},& \widetilde{D_0},& \widetilde{\Delta},&
\widetilde{\Gamma}.
\end{array}$$
A $\Z$-basis for $T_{Km(C\times D)}$ is
$\omega_{1,4}+\omega_{2,3}$, $\omega_{1,3}-\omega_{2,4}.$ and $T_{Km(C\times D)}\simeq \langle 4\rangle^2$.\\
Let $\Dh_4=\langle \sigma_4, \varsigma_2\rangle$. The lattice
$NS(Km(C\times D))^{\Dh_4}$ is generated by the classes
$$\begin{array}{ccccc}
\bar{K}_{W_1} + \bar{K}_{W_2},& \omega_{3,4},& \widetilde{D_0},&
K_{0000} + K_{0001} + K_{0010} + K_{0011},& \hat{K}.
\end{array}$$
The lattice $\Omega_{\Dh_4}:=(H^2(Km(C\times
D))^{\Dh_4})^{\perp}=(NS(Km(C\times D))^{\Dh_4})^{\perp}$ is
generated by the classes
$$
\begin{array}{ccccc}
\widetilde{\Gamma}-2K_{0010} - \widetilde{C_0} - \widetilde{D_0},&
\widetilde{\Delta}-2K_{0010} - \widetilde{C_0} - \widetilde{D_0}
,&
K_{0001} - K_{0010},\\
K_{0000} - K_{0010},& \hat{K} - 2K_{0010} - 4K_{0100} -
2K_{1100},&
\bar{K}_{W_1}-2K_{0010} - 2K_{0100} ,\\
\bar{K}_{W_2}-2K_{0010} - 2K_{0100},&\bar{K}_{W_3} -K_{0010} -
2K_{0100} - K_{1100}  ,&
\bar{K}_{W_4}-K_{0010} - 2K_{0100} - K_{1100},\\
K_{1000} -K_{0100} ,&
K_{1010}  -K_{0100} ,& K_{0110} -K_{0100} ,\\
K_{1001} -K_{0100},&K_{0101} -K_{0100} ,& K_{0011} -K_{0010}.
\end{array}
$$
The lattice $\Omega_{\Dh_4}$ has rank 15, its discriminant is
$-1024$
and its discriminant group is $(\Z/4\Z)^5$.\\
The lattice $H^2(Km(C\times D),\Z)^{\Dh_4}$ is isometric to
$NS(Km(C\times D))^{\Dh_4}\oplus T_{Km(C\times D)}.$
\end{prop}
\bprf The results on $NS(Km(C\times D))$ and $T_{Km(C\times D)}$
are trivial consequences of Remark \ref{rem: divisible classes in
Km(CD)}, of the description of the Kummer lattice given, for
example, in \cite[Chapter VIII, Section 5]{bpv}, \cite[Appendix
5]{Pjateckii Safarevic torelli theorem K3} and more in general of
the known description of the classes generating the second
cohomology group of a Kummer surface (see for example,
\cite{alinikulin}).\\
A direct computation on the coordinates $(x,y,\tau)$ shows that
the automorphism $\varphi_C^3\times \varphi_D$ induces the
automorphism $\sigma_4$ on $Km(C\times D)$. The automorphism
$\varphi_C^3\times \varphi_D$ acts on $C\times D$ as
$(z_1,z_2)\mapsto (-iz_1,iz_2+1)$. So the action of $\sigma_4$ on
the curves of the lattice $K$ (induced by the one of
$\varphi_C^2\times \varphi_D$ on $(C\times D)[2]$) is
$$\begin{array}{ccccccc||ccccccc}
K_{0,0,0,0}\!\!\!\!&\ra \!\!\!\!&K_{0,0,0,1}\!\!\!\!&\ra
\!\!\!\!&K_{0,0,1,1}\!\!\!\!&\ra \!\!\!\!&K_{0,0,1,0}\ \!\!\!\!&
K_{0,1,0,0}\!\!\!\!&\ra \!\!\!\!&K_{1,0,0,1}\!\!\!\!&\ra \!\!\!\!&K_{0,1,1,1}\!\!\!\!&\ra \!\!\!\!&K_{1,0,1,0}\\
K_{1,0,0,0}\!\!\!\!&\ra \!\!\!\!&K_{0,1,0,1}\!\!\!\!&\ra
\!\!\!\!&K_{1,0,1,1}\!\!\!\!&\ra \!\!\!\!&K_{0,1,1,0}\
\!\!\!\!&K_{1,1,0,0}\!\!\!\!&\ra \!\!\!\!&K_{1,1,0,1}\!\!\!\!&\ra
\!\!\!\!&K_{1,1,1,1}\!\!\!\!&\ra
\!\!\!\!&K_{1,1,1,0}.\end{array}$$ Since
$$\sigma_4(\omega_{1,2})=\omega_{1,2},\
\sigma_4(\omega_{1,3})=-\omega_{2,4},\
\sigma_4(\omega_{1,4})=\omega_{2,3},\
\sigma_4(\omega_{2,3})=\omega_{1,4},\
\sigma_4(\omega_{2,4})=-\omega_{1,3},\
\sigma_4(\omega_{3,4})=\omega_{3,4},$$ the classes
$\pi_*(\beta^*[C])$, $\pi_*(\beta^*[D])$ are fixed and
$$\begin{array}{ll}\sigma_4\left(\pi_*(\beta^*[\Delta])\right)=&-\pi_*(\beta^*[\Delta])+2\pi_*(\beta^*[C])+2\pi_*(\beta^*[D]),\\
\sigma_4\left(\pi_*(\beta^*[\Gamma])\right)=&-\pi_*(\beta^*[\Gamma])+2\pi_*(\beta^*[C])+2\pi_*(\beta^*[D]).\end{array}$$
Let $\psi_D$ be the automorphism, of the curve $D$,
$\psi_D:(\tau,\nu)\ra (1/\tau,i\nu/\tau^2)$. The automorphism
$\varphi_C^3\times \psi_D$ induces on $Km(C\times D)$ the
automorphism $\varsigma_2$. We observe that $\varphi_C^3$ and
$\psi_D$ have order four respectively on $C$ and $D$. However
$\varphi_C^3\times \psi_D$ induces an automorphisms of order two
on $Km(C\times D)$, indeed $(\varphi_C^3\times
\psi_D)^2=\iota_C\times \iota_D$ which is the identity on
$Km(C\times D)$.\\ The automorphism $\psi_D$ of $D$ has order four
and fixes the points $0_D=(1,0)$ and $(-1,0)$. Hence it
corresponds to $z_2\mapsto iz_2$. So the automorphism
$\varsigma_2$ is induced by the automorphism $(z_1,z_2)\mapsto
(-iz_1,iz_2)$ on $C\times D$. The action of $\varsigma_2$ on the
curves of the Kummer lattice is
$$\begin{array}{ccccccccccc}
K_{0,0,0,1}&\leftrightarrow &K_{0,0,1,0};& &K_{0,1,0,0}&\leftrightarrow &K_{1,0,0,0};&  &K_{0,1,0,1}&\leftrightarrow &K_{1,0,1,0};\\
K_{0,1,1,0}&\leftrightarrow
&K_{1,0,0,1};&&K_{0,1,1,1}&\leftrightarrow &K_{1,0,1,1};&
&K_{1,1,0,1}&\leftrightarrow &K_{1,1,1,0}\end{array}$$ and the
other classes of the Kummer lattice are fixed. The action of
$\varsigma_2$ on the classes $\pi_*(\beta^*[C])$,
$\pi_*(\beta^*[D])$, $\pi_*(\beta^*[\Delta])$,
$\pi_*(\beta^*[\Gamma])$ is the same as the action of
$\sigma_4$.\\
Since we described the action of $\sigma_4$ and $\varsigma_2$, it
is trivial to compute $\Omega_{\Dh_4}$ and its orthogonal in
$H^2(Km(C\times D),\Z)$.\erem

\begin{prop}\label{prop: NS and MW elliptic with dihedral}
Let $X$ be the generic elliptic K3 surface admitting the
automorphisms $\sigma_4$ and $\varsigma_2$, i.e. $X$ is the
generic K3 surface with equation \eqref{formula: equation elliptic
with dihedral}. Then $NS(X)\simeq U\oplus \Omega_{\Dh_4}$ and
$MW(X)\simeq \Omega_{\Dh_4}$.
\end{prop}
\bprf In \cite{Whitcher} is proved that the action of a finite
group of symplectic automorphisms on the second cohomology group
of the K3 surface (and hence on $\Lambda_{K3}$) is essentially
unique (cf. \cite[Definition 4.6]{Nikulin symplectic}). Since the
action of $G$ on $H^2(X,\Z)$ does not depend on $X$ up to
isometry, the lattice $(H^2(X,\Z)^G)^{\perp}=:\Omega_G$ does not
depend on $X$. In particular this implies that if $G$ is a finite
group of symplectic automorphisms of a K3 surface $S$, the lattice
$\Omega_G$ is primitively embedded in $NS(S)$, indeed since $G$ is
symplectic $T_S\hookrightarrow H^2(S,\Z)^G$ and hence
$NS(S)\hookleftarrow(H^2(S,\Z)^G)^{\perp}$. In Proposition
\ref{prop: action of symplectic dihedral} we computed the lattice
$\Omega_{\Dh_4}$. It follows that $\Omega_{\Dh_4}\hookrightarrow
NS(X)$. The proposition follows as in the proof of Proposition
\ref{prop: symplectic base change}.\erem
\begin{rem}{\rm Proposition \ref{prop: NS and MW elliptic with
dihedral} gives the N\'eron--Severi group and the Mordell--Weil
lattice of the elliptic K3 surface admitting 15 independent
sections described in \cite[Theorem 1.1]{Kloosterman} (cf. Reamrk
\ref{rem: kloosterman}).\erem}\end{rem}

The paper \cite{EE8- lattices and dihedral} classifies lattices
which have dihedral groups $\Dh_n$ (for $n=2,3,4,5,6$) in the
group of the isometries and which have the properties:
\begin{itemize}\item the lattices are rootless (i.e. there are no elements of length $2$),
\item they are the sum of two copies of $E_8(2)$, \item there are
involutions in $\Dh_n$ acting as minus the identity on each copy
of $E_8(2)$.\end{itemize} In \cite[Section 5.2]{EE8- lattices and
dihedral} it is proved that there are three lattices satisfying
these properties and admitting $\Dh_4$ as group of isometries,
they are called $DIH_8(15)$, $DIH_8(16,DD4)$, $DIH_8(16,0)$.

\begin{rem}\label{rem: omegadh4 is isomteric to DIH8(16)}{\rm The lattice $DIH_8(15)$ described in \cite{EE8- lattices and
dihedral} is isometric to the lattice $\Omega_{\Dh_4}(-1)$.\\
\binf $\Omega_{\Dh_4}$ has no vectors of length $-2$ (this is a
property of all the lattices $\Omega_{G}$, cf. \cite[Lemma
4.2]{Nikulin symplectic}, \cite{Whitcher}). The two involutions
$\varsigma_2$ and $\varsigma_2 \sigma_4$ generate the group
$\Dh_4$ as group of isometries on $\Omega_{\Dh_4}$. Moreover they
act as minus identity on two copies of $E_8(-2)$ (because they are
both symplectic involutions, so they act as minus the identity on
$\Omega_{\Z/2\Z}$ and $\Omega_{\Z/2\Z}\simeq E_8(-2)$, cf.
\cite[Section 1.3]{bert Nikulin involutions}) and the sum of these
two copies of $E_8(-2)$ is $\Omega_{\Dh_4}$ (this follows from the
explicit computation on the action of $\varsigma_2$ and
$\varsigma_2\sigma_4$ on $\Omega_{\Dh_4}$), hence
$\Omega_{\Dh_4}(-1)$ is isometric to one of the lattices
$DIH_8(15)$, $DIH_8(16,DD4)$, $DIH_8(16,0)$. Since
$\Omega_{\Dh_4}$ has rank 15 and the lattices $DIH_8(15)$,
$DIH_8(16,DD4)$, $DIH_8(16,0)$ have rank 15, 16, 16 respectively,
$\Omega_{\Dh_4}(-1)\simeq DIH_8(15)$.\erem}\end{rem}

\section{Dihedral groups of automorphisms on elliptic fibrations with
torsion.}\label{section: dihedral groups of automorphisms on
elliptic fibration with torsion} In \cite{dihedral 5} we proved
that if a K3 surface admits $G=\Z/5\Z$ as group of symplectic
automorphisms, then the same surface admits $\Dh_5$ as group of
symplectic automorphisms. Moreover it is proved that
$\Omega_{\Z/5\Z}$ is isometric to $\Omega_{\Dh_5}$. A similar
phenomenon happens also for other groups of symplectic
automorphisms. Here we will consider these other cases and we
reconsider the situation described in \cite{dihedral 5}.\\
We consider the following groups with their presentation:
$$\begin{array}{c}\Dh_5=\langle c,d| c^5=1,\ d^2=1, cd=dc^{-1}  \rangle, \\
\Dh_6=\langle c,d| c^6=1,\ d^2=1, cd=dc^{-1}\rangle\\
\Z/2\Z\times \Dh_4=\langle c,d,e| c^2=1,\ d^2=1, e^4=1,\ ce=ec,\
cd=dc,\ ed=de^{-1}\rangle\\ J:=\langle c,d,e| c^2=1,\ d^3=1,
e^3=1,\ cd=d^{-1}c,\ ce=e^{-1}c\rangle.\end{array}$$  Let $G$ be
one of the following four groups: $\Z/n\Z$, $n=5,6$, $\Z/2\Z\times
\Z/4\Z$, $(\Z/3\Z)^2$. In \cite{symplectic not prime} the equation
of the generic elliptic K3 surface $X_G$ admitting $G$ as torsion
part of Mordell--Weil group is given. For these groups the
equation of elliptic fibration $\mathcal{E}_G$ on $X_G$ depends on
the choice of two polynomials of degree 2. More precisely the
equation of $\mathcal{E}_G$ is of type
$$y^2=x^3+C(a,b)x^2+A(a,b)x+B(a,b)$$
where $a(\tau)=\widetilde{a}(\tau,1)$,
$b(\tau)=\widetilde{b}(\tau,1)$, $\widetilde{a}(\tau,\sigma)$,
$\widetilde{b}(\tau,\sigma)$ are two homogenous polynomials of
degree two in $\tau$, $\sigma$, and $A,B,C$ are homogenous
polynomials in $a,b$ of degree 4,6,2 respectively (cf.
\cite[Tables 1 and 2]{symplectic not prime}).

\begin{prop}\label{prop: XG admit also a dhiedral group of automorphisms} Let
$X_G$ be the K3 surface admitting the generic elliptic fibration
$\mathcal{E}_G$ such that $MW(\mathcal{E}_G)=G$. Let $G$ be one of
the four groups $\Z/n\Z$, $n=5,6$, $\Z/2\Z\times \Z/4\Z$,
$(\Z/3\Z)^2$ and let $t_1$ be the a generator of $\Z/n\Z$,
$n=5,6$, of $\Z/4\Z$, of a copy of $\Z/3\Z$ respectively. Then
there exists a symplectic involution $\iota$ on $X_G$ and this
involution is such that $\Dh_m=\langle
\sigma_{t_1},\iota \rangle$ where $m$ is the order of $t_1$.\\
Moreover \begin{equation}\label{equation: isometries between
omegaG e altri omega}\Omega_{\Z/n\Z}\simeq \Omega_{\Omega_n}\mbox{
for }\ n=5,6,\ \ \ \ \Omega_{\Z/2\Z\times \Z/4\Z}\simeq
\Omega_{\Z/2\Z\times \Dh_4},\ \ \ \Omega_{\Z/3\Z\times
\Z/3\Z}\simeq \Omega_{J}.\end{equation}
\end{prop}
\bprf First we will prove that the elliptic fibration
$\mathcal{E}_G$ admits a symplectic involution. To construct this
involution, we consider an automorphism on $\mathbb{P}^1$.\\
\textit{Step 1: construction of an involution on $\mathbb{P}^1$.}
Let $\widetilde{a}(\tau,\sigma)$, $\widetilde{b}(\tau,\sigma)$ be
two homogeneous polynomials of degree two. Let $\alpha_1$ and
$\alpha_2$ be the roots of $\widetilde{a}$ and $\beta_1$,
$\beta_2$ the ones of $\widetilde{b}$. There is an involution
$\theta$ of $\mathbb{P}^1$ such that $\theta(\alpha_1)=\alpha_2$
and $\theta(\beta_1)=\beta_2$. Indeed up to an isomorphism of
$\mathbb{P}^1$ one can always suppose that
$\widetilde{a}(\tau,\sigma)=\tau\sigma$ and
$\widetilde{b}(\tau,\sigma)=\lambda(\tau-\sigma)(\tau-\mu\sigma)$.
Hence the involution $\theta(\tau,\sigma)=(\mu \sigma,\tau)$ has
the required property. We note that on the affine coordinate
$\tau$ the action of $\theta$ is $\tau\mapsto \mu/\tau$.\\
\textit{Step 2: construction of a symplectic involution on $X_G$.}
The automorphism
$$\vartheta:(x,y,\tau)\mapsto
(\frac{\mu^2x}{\tau^4},-\frac{\mu^3y}{\tau^6},\frac{\mu}{\tau})=(\frac{\mu^2x}{\tau^4},-\frac{\mu^3y}{\tau^6},\theta(\tau))$$
is a symplectic automorphism of $X_G$. Indeed, $\mathcal{E}_G$ has
an equation of type:
$$y^2=x^3+C(a(\tau),b(\tau))x^2+A(a(\tau),b(\tau))x+B(a(\tau),b(\tau)).$$
Since $a(\theta(\tau))=\frac{\mu}{\tau^2}a(\tau)$ and
$b(\theta(\tau))=\frac{\mu}{\tau^2}b(\tau)$ it is clear that the
equation is invariant under $\vartheta$. Moreover, this
automorphism is symplectic, in fact $\vartheta^*(dx/y\wedge
d\tau)=dx/y\wedge d\tau$.\\
\textit{Step 3: description of the group $\langle
\sigma_{t_1},\iota\rangle$} Let $t_1$ be a section of order
$5,6,4,3$ if $G$ is $\Z/5\Z$, $\Z/6\Z$, $\Z/2\Z\times\Z/4\Z$,
$(\Z/3\Z)^2$ respectively. Let $\iota=\sigma_{t_1}\circ\vartheta$,
then $\iota$ is a symplectic automorphism on $X_G$, because it is
composition of two symplectic automorphisms. Let $F_s$ be the
fiber over the point $s\in \mathbb{P}^1$. The fibers $F_{\tau}$
and $F_{\iota_{\mathbb{P}^1}(\tau)}$ are isomorphic and we chose
an isomorphism such that
$t_i(\tau)=t_i(\iota_{\mathbb{P}^1}(\tau))$ for each section
$t_i$, hence to each point $P_{\tau}\in F_{\tau}$ corresponds
uniquely a point $P_{\iota_{\mathbb{P}^1}(\tau)}\in
F_{\iota_{\mathbb{P}^1}(\tau)}$. The automorphism $\iota$ acts in
the following way $P_{\tau}\mapsto
-P_{\iota_{\mathbb{P}^1}(\tau)}+T_{\iota_{\mathbb{P}^1}(\tau)}$
where '$-$' is the operation on the elliptic curve and $T_{s}$ is
the intersection of the fiber $F_s$ with the section $t_1$. Its
clear that $\iota$ is an involution and its action on the section
is $\iota(t_i)=t_{-i+1}$. This involution, together with
$\sigma_{t_1}$, generates the dihedral group $\Dh_n$ if $t_1$ is
an $n$-torsion
section.\\
Moreover, if $G=\Z/2\Z/\times\Z/4\Z$ and $r$ is the 2-torsion
section such that $r\notin\langle t_1\rangle$, then
$\iota\sigma_r=\sigma_r\iota$ and hence
$\langle\sigma_r,\sigma_t,\iota\rangle=\Z/2\Z\times \Dh_4$.\\
If $G=(\Z/3\Z)$ and $r$ is a 3-torsion section such that
$r\notin\langle t_1\rangle$, then
$\iota\sigma_r=\sigma_r^{-1}\iota$, hence
$\langle\sigma_r,\iota\rangle\simeq\langle\sigma_t,\iota\rangle=\Dh_3$
and $\langle\sigma_r,\sigma_t,\iota\rangle\simeq J$.\\
\textit{Step 4: proof of \eqref{equation: isometries between
omegaG e altri omega}.} The automorphism $\iota$ fixes the class
of the fiber on $\mathcal{E}_G$ and the sum of all the sections
(it clearly sends sections to sections). Since $\iota$ does not
fix the points of $\mathbb{P}^1$ corresponding to singular fibers
(this follows from the computation of the discriminant of the
elliptic fibration $\mathcal{E}_G$), there are no components of
the reducible fibers which are fixed by $\iota$. Hence
$$NS(X)^{\iota}\otimes \Q=\langle F,\sum_{t\in
MW(\mathcal{E}_G)}t\rangle\otimes \Q\ \ \mbox{ and so }\ \
(NS(X)^{\iota})^{\perp}=\langle F,\sum_{t\in
MW(\mathcal{E}_G)}t\rangle^{\perp}.$$ We recall that also
$NS(X)^{G}\otimes \Q=\langle F, \sum_{t\in
MW(\mathcal{E}_G)}t\rangle\otimes \Q$ (cf. \cite{symplectic
prime}, \cite{symplectic not prime}) and hence $$\Omega_G\simeq
(NS(X)^{G})^{\perp}\simeq \langle F, \sum_{t\in
MW(\mathcal{E}_G)}t\rangle^\perp.$$ Let $H=\langle \iota,
G\rangle$. Since $NS(X)^G\otimes \Q=NS(X)^{\iota}\otimes \Q$, we
have $NS(X)^H\otimes \Q=NS(X)^G\otimes \Q$ and so $$\Omega_H\simeq
(NS(X)^H)^{\perp}\simeq(NS(X)^G)^{\perp}\simeq \Omega_G.$$\erem

\begin{rem}{\rm The group $J$ is isomorphic to the group
$\Ag_{3,3}$ described in \cite{mukai}. Indeed $J$ is a finite
group of symplectic automorphisms on a K3 surface. These groups
are listed by Xiao, \cite[Table 2]{Xiao}. The order of $J$ is 18
and there are only two groups in the Xiao's list with order 18.
One of them is $\Z/3\Z\times\Dh_3 $ but in $J$ there are no
elements of order 3 commuting with an involution, hence $H$ has to
be isomorphic to the other one, which is
$\Ag_{3,3}.$\erem}\end{rem}

In \cite{dihedral 5} the case of the elliptic fibration
$\mathcal{E}_{\Z/5\Z}$ is analyzed. In this situation it is
possible to show that $E_8(-1)\oplus E_8(-1)\subset
NS(X_{\Z/5\Z})$. Hence one can define an involution which switches
these two copies of $E_8(-1)$ (by Definition this is a
Morrison--Nikulin involution, cf. \cite[Section 4.4]{bert Nikulin
involutions}). In the case $G=\Z/6\Z$ the situation is very
similar.\\ We will need the following trivial lemma and so we
recall that given a lattice $T$, the length $l(A_T)$ is is the
minimum number of generator of $A_T:=T^{\vee}/T$.

\begin{lemma}\label{lemma: condition on T=L(n)}
Let $n\in \N$ and $(T,b_T$) be a lattice such that: {\it i)}
rk$T=l(A_T)$;\\ {\it ii)} $T^{\vee}/T\simeq \bigoplus_i \Z/nd_i\Z$
and $\beta_i$
generate $\Z/nd_i\Z$ in $T^{\vee}/T$.\\
Let $L$ be the free $\Z$-module $L=\{x\in T\otimes \Q\mbox{ such
that }nx\in T\}$ and $b_L=n(b_{T\otimes\Q})_{|L}$ a bilinear form
defined on $L$.
Then $(L,b_L)$ is a lattice and $T\simeq L(n)$ if and only if $b_{T^\otimes \Q}(d_i\beta_i,d_j\beta_j)\in \frac{1}{n}\Z$.\\
Moreover $L$ is an even lattice if and only if
$b_{T\otimes\Q}(d_i\beta_i,d_i\beta_i)\in \frac{2}{n}\Z$.
\end{lemma}

\begin{prop}\label{prop: quotient of E by iota and Dhn} Let $G$ be
one of the following four groups: $\Z/n\Z$, $n=5,6$, $\Z/2\Z\times
\Z/4\Z$, $(\Z/3\Z)^2$. Let $\mathcal{E}_{G}$ be the generic
elliptic fibration such that $MW(\mathcal{E}_{G})=G$. Then: the
involution $\iota$ (defined in proof of Proposition \ref{prop: XG
admit also a dhiedral group of automorphisms}) is a
Morrison--Nikulin involution for $G=\Z/n\Z$, $n=5,6$ and
$G=(\Z/3\Z)^2$ and it is {\bf not} a Morrison--Nikulin involution
if $G=\Z/2\Z\times\Z/4\Z$.
\end{prop}
\bprf {\it The case $G=\Z/n\Z$, $n=5,6$.} We consider only the
case $n=6$ because the other is very similar and in part it is
analyzed in \cite{dihedral 5}. The equation of
$\mathcal{E}_{\Z/6\Z}$ is given \cite[Table 1]{symplectic not
prime}. We recall that the reducible fibers of the fibration
$\mathcal{E}_{\Z/6\Z}$ are $2 I_6+2I_3+2I_2+2I_1$ (we will number
the reducible fibers assuming that the first two reducible fibers
are of type $I_6$, the third and the fourth are of type $I_3$, the
fifth and the sixth are of type $I_2$). Let $t_1$ be a 6-torsion
section, then $t_1\cdot C_{1}^{(j)}=1$ for all $j=1,\ldots, 6$ and
$t_1\cdot C_i^{(j)}=0$ if $i\neq 1$. Here we give two orthogonal
copies of $E_8(-1)$ which are contained in $NS(X_{\Z/6\Z})$:
\begin{eqnarray*}
\begin{array}{ccccccccccccc}
C_2^{(3)}&-&C_0^{(3)}&-&s&-&C_0^{(1)}&-&C_5^{(1)}&-&C_4^{(1)}&-&C_3^{(1)}\\
      & &      & &\mid & &   & &   & &   & &\\
      & &      & & C_0^{(5)}   & &   & &   & &   & &\\
\\
C_2^{(4)}&-&C_{1}^{(4)}&-&t_1&-&C_{1}^{(2)}&-&C_2^{(2)}&-&C_{3}^{(2)}&-&C_{4}^{(2)}\\
      & &      & &\mid & &   & &   & &   & &\\
      & &      & & C_{1}^{(6)}   & &   & &   & &   & &
\end{array}
\end{eqnarray*}
There exists a Morrison--Nikulin involution, $\iota_{M}$, of the
surface which switches them (cf. \cite[Proposition
5.7]{morrison}). This involution switches the singular fibers of
the same type and acts on $MW(\mathcal{E}_{\Z/6\Z})$ sending
$s\leftrightarrow t_1$, $t_2\leftrightarrow t_5$,
$t_3\leftrightarrow t_5$ (this is deduced by the fact that
$t_i\cdot C_j^{(h)}=\iota_M(t_i)\cdot \iota_M(C_j^{(h)})$). We
observe that the action of $\iota_M$ on the sections is
$t_i\leftrightarrow t_{-i+1}$. Hence $\iota_M$ is the involution
defined in proof of Proposition \ref{prop: XG admit also a
dhiedral group of automorphisms}. In particular $\Dh_6=\langle
\iota_M,\sigma_6\rangle$.\\
{\it The case $G=(\Z/3\Z)^2$.} The elliptic fibration
$\mathcal{E}_{(\Z/3\Z)^2}$ has 8 fibers of type $I_3$ (cf.
\cite[Table 2]{symplectic not prime}). The transcendental lattice
of the surface $X_{(\Z/3\Z)^2}$ is isometric to $U(3)\oplus U(3)$.
The automorphisms $\iota$ switches pairwise of fibers of type
$I_3$. On the quotient elliptic fibration this gives 4 fibers of
type $I_3$. Moreover the automorphism $\iota$ fixes two smooth
fibers, and on each of these fibers it has four fixed points. On
the quotient elliptic fibration $\mathcal{E}/\iota$ this gives two
fibers of type $I_0^*$. So the trivial lattice of the elliptic
fibration $\mathcal{E}/\iota$ is $U\oplus A_2^{\oplus 4}\oplus
D_4^{\oplus 2}$ and there are no sections of infinite order,
because
$rk(Tr_{\mathcal{E}/\iota})=rk(NS(\widetilde{X_{(\Z/3\Z)^2}/\iota}))$.
The discriminant of the trivial lattice is $3^42^4$. By
\cite{shimada} the elliptic fibration $\mathcal{E}/\iota$ is such
that $tors(MW(\mathcal{E}_{(\Z/3\Z)^2/\iota}))=\{1\}$. This
implies that $d(NS(\widetilde{X_{(\Z/3\Z)^2}/\iota}))=3^42^4=6^4$.
Moreover it is easy to find the discriminant form of
$NS(\widetilde{X_{(\Z/3\Z)^2}/\iota})$. The discriminant form of
the transcendental lattice is the opposite of the one of the
N\'eron--Severi group. Hence one can directly check that the
hypothesis of Lemma \ref{lemma: condition on T=L(n)} with $n=6$
are satisfied by $T_{\widetilde{X_{(\Z/3\Z)^2}/\iota}}$, and so
there exist a lattice $L$ such that
$T_{\widetilde{X_{(\Z/3\Z)^2}/\iota}}=L(6)$. The discriminant form
of $L$, computed by the one of
$T_{\widetilde{X_{(\Z/3\Z)^2}/\iota}}$, is trivial, hence $L$ is a
unimodular lattices with signature $(2,2)$. This implies that
$L\simeq U\oplus U$ and so
$T_{\widetilde{X_{(\Z/3\Z)^2}/\iota}}\simeq U(6)\oplus U(6)$.
Since $U(6)\oplus U(6)\simeq (T_{X_{(\Z/3\Z)^2}})(2)$, the
involution $\iota$ is a Morrison--Nikulin involution (cf.
\cite[Theorems 5.7, 6.3]{morrison}).\\
{\it The case $G=\Z/2\Z\times \Z/4\Z$.} The elliptic fibration
$\mathcal{E}_{\Z/2\Z\times\Z/4\Z}$ has 4 fibers of type $I_4$ and
4 fibers of type $I_2$ (cf. \cite[Table 2]{symplectic not prime}).
The automorphisms $\iota$ switches two fibers of type $I_4$, two
other fibers of type $I_4$, two fibers of type $I_2$ and the
remaining two fibers of type $I_2$. Moreover it fixes two smooth
fibers and in particular four distinct points on each of these
fibers. As in case $G=(\Z/3\Z)^2$ this gives the quotient elliptic
fibration. The trivial lattice of the elliptic fibration
$\mathcal{E}/\iota$ is $U\oplus A_1^{\oplus 2}\oplus A_3^{\oplus
2}\oplus D_4^{\oplus 2}$ and there are no sections of infinite
order. The discriminant of the trivial lattice is $2^{10}$. By
\cite{shimada} the elliptic fibration $\mathcal{E}/\iota$ is such
that
$tors(MW(\mathcal{E}_{\Z/2\Z\times\Z/4\Z}/\iota))=(\Z/2\Z)^2$.
This implies that $d(NS(\widetilde{X/\iota}))=2^{10}/2^4=2^6$. We
observe that also $d(NS(X_{\Z/2\times\Z/4\Z}))=2^6$. In particular
this implies that $\iota$ is not a Morrison--Nikulin involution,
indeed, if $Y=\widetilde{X/i}$ with a Morrison--Nikulin involution
$i$, then $T_Y=T_X(2)$ (cf. \cite[Theorems 5.7, 6.3]{morrison}) so
in particular $d(NS(X))=2^{22-\rho(X)}d(NS(Y))$, where $\rho(X)$
is the Picard number of $X$.\erem

\begin{prop}\label{prop: quotient of E by iota and Dhn} Let $G$, $\iota$ and $\mathcal{E}_{G}$ be as in Proposition \ref{prop: XG admit also a dhiedral group of automorphisms}. Let $K=\langle
G,\iota\rangle$. Then:
\begin{itemize} \item[i)] the
desingularization, $\widetilde{X_{G}/K}$, of the quotient
$X_{G}/K$ is such that $NS(\widetilde{X_{G}/K})\simeq
NS(\widetilde{X_G/\iota})$; \item[ii)] the desingularization of
the quotient $\widetilde{X_G/K}$ is such that
$NS(\widetilde{X_G/K})$ is an overlattice of index 4 of
$U(2)\oplus M_{K}$.\end{itemize} If $G=\Z/n\Z$, $n=5,6$, then
$T_{\widetilde{X_{G}/\iota}}\simeq T_{\widetilde{X_{G}/K}}\simeq
U(2)\oplus U(2n)$, if $G=(\Z/3\Z)^2$, then
$T_{\widetilde{X_{G}/\iota}}\simeq T_{\widetilde{X_{G}/K}}\simeq
U(6)\oplus U(6).$
\end{prop}
\bprf Here we consider the case $n=6$. The others are very
similar.\\
{\it Proof of i).} Let $\mathcal{E}/\iota$ be the quotient
elliptic fibration. The involution $\iota$ switches the two fibers
of type $I_6$, the two fibers of type $I_3$, the two fibers of
type $I_2$. On $\mathcal{E}/\iota$ these reducible fibers
correspond to one fiber of type $I_6$, one of type $I_3$, one of
type $I_2$. On the fibration $\mathcal{E}$, $\iota$ fixes two
smooth fibers. On these fibers $\iota$ fixes 4 points (the points
$P$ such that $2P=T$ where $T$ is the intersection of the fiber
with the section $t_1$). Hence these two fibers corresponds on
$\mathcal{E}/\iota$ to 2 fibers of type $I_0^*$. So the trivial
lattice of the fibration $\mathcal{E}/\iota$ is $U\oplus A_5\oplus
A_2\oplus A_1\oplus D_4^{\oplus 2}$. Choosing a numbering of the
fibers and of the irreducible components ($D_i^{(j)}$) of the
fibers we assume that the first fiber is of type $I_6$, the second
of type $I_3$, the third of type $I_2$, the fourth and the fifth
of type $I_0^*$ and $F\simeq
D_0^{(4)}+D_1^{(4)}+2D_2^{(4)}+D_3^{(4)}+D_4^{(4)}\simeq
D_0^{(5)}+D_1^{(5)}+2D_2^{(5)}+D_3^{(5)}+D_4^{(5)}$ where $F$ is
the class of the fiber and $'\simeq'$ is the linear equivalence of
divisors.\\ The elliptic fibrations with trivial lattice $U\oplus
A_5\oplus A_2\oplus A_1\oplus D_4^{\oplus 2}$ admit a 2-torsion
section $r$ (cf. \cite{shimada}) and up to a choice of numbering
of the components of the reducible fibers we can assume, by the
height formula, that the intersection numbers are:
$$r\cdot D_3^{(1)}=r\cdot D_0^{(2)}=r\cdot D_1^{(3)}=r\cdot D_1^{(4)}=r\cdot D_1^{(5)}=1$$
and the other intersections are zero.\\
Now let us consider the quotient $X/\Dh_6$. Under $\Dh_6$, the
orbit of the curve $C_0^{(1)}$ consists of all the curves
$C_i^{(j)}$, $i=0,\ldots,5$, $j=1,2$, so it is made up of all the
components of the two fibers of type $I_6$. Hence the image of
these two fibers on $X/\Dh_6$ is a nodal curve, where the node is
the image of the singular points of the two fibers of type $I_6$.
Analogously, the image of the two fibers of type $I_3$ is a nodal
curve. The node is the image of the singular points of these
fibers, and it is a singularity of type $A_1$ for the surface
$X/\Dh_6$. Indeed on $X$ the singular points of the fibers of type
$I_3$ are fixed by $\sigma_{t_1}^2$ and hence their stabilizer is
$\Z/2\Z$, this implies that the corresponding point on the
quotient is a singularity of type $A_1$. With the same argument,
one notices that the two fibers of type $I_2$ correspond on the
quotient to a nodal curve, whose node is a singular point for
$X/\Dh_6$ of type $A_2$ and the two fibers of type $I_1$
correspond on the quotient to a nodal curve, whose node is a
singular point for $X/\Dh_6$ of type $A_5$. Each of the two smooth
fibers fixed by $\iota$ correspond on the quotient to a rational
curves with 4 singularities of type $A_1$. Blowing up these
singularities we find the quotient elliptic fibration
$\mathcal{E}/\Dh_6$. Its singular fibers
are$I_1+I_2+I_3+I_6+2I_0^*$. The trivial lattice of this fibration
is $U\oplus A_1\oplus A_2\oplus A_5\oplus D_4^{\oplus 2}$. It has
a zero section (the image of the sections of $\mathcal{E}$), it
has no sections of infinite order (because it has Picard number
equal to 18, as the one of $X$ and its trivial lattice has rank
18) and has a 2-torsion section (\cite{shimada}).The singular
fibers of the $\mathcal{E}/\Dh_6$ are the same as the ones of
$\mathcal{E}/\iota$ and
$MW(\mathcal{E}/\Dh_6)=MW(\mathcal{E}/\iota)$. Moreover up to a
choice of the numbering on the components of the reducible fibers
one can assume that the intersection property of the 2-torsion
section on $\mathcal{E}/\Dh_6$ are the same as the one of the
2-torsion section of $\mathcal{E}/\iota$. The N\'eron--Severi
group of an elliptic surface is completed determined by its
singular fibers and by its sections and hence
$NS(\widetilde{X/\iota})\simeq NS(\widetilde{X/\Dh_6})$.\\
{\it Proof of ii).} Let us consider again the quotient
$\widetilde{X/\Dh_6}$. We will call the components of its
reducible fibers $B_i^{(j)}$ and we number the fibers assuming
that the first is of type $I_6$, the second of type $I_2$, the
third of type $I_3$ and the fourth and fifth of type $I_0^*$. By
the construction described above the components which do not
intersect the zero section of the fibers of type $I_6$, $I_3$,
$I_2$ arise from the desingularization of singular points of type
$A_5$, $A_2$, $A_1$ respectively on the quotient $X/\Dh_6$. So the
components $B_i^{(1)}$, $i=1,\dots ,5$, $B_j^{(2)}$, $j=2,3$ and
$B_1^{(3)}$ are curves in $M_{\Dh_6}$. The fibers of type $I_0^*$
on $\widetilde{X/\Dh_6}$ are obtained as desingularization of a
rational curve which has 4 singularity of type $A_1$, hence the
components $B_i^{(j)}$ of the fourth and fifth fibers which
resolve the singularity of type $A_1$ on the quotient $X/\Dh_6$
are $B_h^{(k)}$, $h=0,1,3,4$, $k=4,5$ and hence they are curves in
$M_{\Dh_6}$. The curves $B_i^{(1)}$, $i=1,\dots ,5$, $B_j^{(2)}$,
$j=2,3$, $B_1^{(3)}$, $B_h^{(k)}$, $h=0,1,3,4$, $k=4,5$ generates
the lattice $A_5\oplus A_2\oplus A_1^{\oplus 9}$ which is a
sublattice of $M_{\Dh_6}$. The lattice $M_{\Dh_6}$ is described in
\cite[Table 2, Lemma 2]{Xiao} as overlattice of $A_5\oplus
A_2\oplus A_1^{\oplus 9}$ obtained by adding two 2-divisible
classes and hence $A_5\oplus A_2\oplus A_1^{\oplus
9}\hookrightarrow M_{\Dh_6}$ is an inclusion with index 4. Indeed
the classes
$$\begin{array}{ll}v:=\frac{1}{2}(B_0^{(4)}+B_1^{(4)}+B_3^{(4)}+B_4^{(4)}+
B_0^{(5)}+B_1^{(5)}+B_3^{(5)}+B_4^{(5)})&\sim
(F+B_2^{(4)}+B_2^{(5)}),\\
w:=\frac{1}{2}\left(B_1^{(1)}+B_3^{(1)}+B_5^{(1)}+B_1^{(3)}+\sum_{j=4}^5\left(B_3^{(j)}+B_4^{(j)}\right)\right)&\sim
r-2F-s,\end{array}$$ where $r$ is the 2-torsion section of
$\mathcal{E}/\Dh_6$, are in $M_{\Dh_6}\subset NS(\widetilde{X/\Dh_6})$.\\
The orthogonal lattice to $M_{\Dh_6}$ in $NS(\widetilde{X/\Dh_6})$
is made up of the two classes $F$ and
$z:=-2s+B_1^{(4)}+2B_2^{(4)}+B_3^{(4)}+B_4^{(4)}+B_1^{(5)}+2B_2^{(5)}+B_3^{(5)}+B_4^{(5)}.$
They generate a lattice which is isometric to $U(2)$. Hence
$U(2)\oplus M_{\Dh_6}\hookrightarrow NS(\widetilde{X}/\Dh_6)$ with
a finite index. This index is $4=\left(d(U(2)\oplus
M_{\Dh_6})\right)/\left(d(NS(\widetilde{X}/\Dh_6)\right)$ and in
fact the classes
$$\begin{array}{lc}
z+B_1^{(4)}+B_3^{(4)}+B_4^{(4)}+B_1^{(5)}+B_3^{(5)}+B_4^{(5)}&\in
U(2)\oplus M_{\Dh_6}\\
F+B_0^{(4)}+B_1^{(4)}+B_3^{(4)}+B_4^{(4)}&\in U(2)\oplus
M_{\Dh_6}\end{array}$$ are divisible by 2 in
$NS(\widetilde{X}/\Dh_6)$.\\
Let $\nu$ be a Morrison--Nikulin involution on a K3 surface $S$.
Then $T_{\widetilde{S/\nu}}\simeq T_{S}(2)$ (cf. \cite[Theorems
5.7 and 6.3]{morrison})). The transcendental lattice of
$X_{\Z/n\Z}$ is computed in \cite{symplectic prime} and is
isometric to $U\oplus U(n)$. So
$T_{\widetilde{X_{\Z/n\Z}/\iota}}\simeq U(2)\oplus U(2n).$ Since
$NS(\widetilde{X/\iota})\simeq NS(\widetilde{X/K})$, the
discriminant form of $T_{\widetilde{X/\iota}}$ is isometric to the
one of $T_{\widetilde{X/K}}$. The lattice with such a discriminant
form is uniquely determined, up to isometry, indeed by Lemma
\ref{lemma: condition on T=L(n)} it is multiple of a certain
lattice which is uniquely determined by \cite[Corollary
1.13.3]{Nikulin bilinear}. Hence
$T_{\widetilde{X_{\Z/n\Z}/\iota}}\simeq
T_{\widetilde{X_{\Z/n\Z}/\Dh_6}}$. \erem

\begin{rem}\label{rem: omegadh5 is isomteric to DIH10(16)}{\rm As in Remark \ref{rem: omegadh4 is isomteric to DIH8(16)}
one obtains that for $n=5,6$ the lattice $DIH_{2n}(16)$ described
in \cite{EE8- lattices and dihedral} is isometric to the lattice
$\Omega_{\Dh_n}(-1)$, and, by Remark \ref{rem: isometries Shioda
lattices}, they are both isometric to $MW(F^{(n)}_{gen})$,
described in \cite{Shioda F5n in generale}. The isometry
$\Omega_{\Dh_5}(-1)\simeq DIH_{10}(16)$ was proved also in
\cite{dihedral 5}.\erem}\end{rem} In Proposition \ref{prop: NS MW
with abelina group} we exhibit K3 surfaces obtained as quotient of
other K3 surfaces $X$ with an abelian group $G$ of symplectic
automorphisms, such that $NS(X)\simeq U\oplus M_G$. In the case
$G=\Dh_5$, $\Dh_6$, $\Ag_{3,3}$, $\Z/2\Z\times \Dh_4$ this is not
possible, indeed there exists no K3 surfaces with N\'eron--Severi
group isometric to $U\oplus M_{G}$ for one of these group $G$.
\section{Final remarks on finite groups of symplectic
automorphisms.}\label{section: final remarks on finite groups of
symplectic automorphisms}

In Section \ref{section: dihedral groups of automorphisms on
elliptic fibration with torsion} we analyzed the following
property of a pair of finite groups $(H,G)$:
\begin{assumeP}\label{property P}
$\bullet$ $G$ acts symplectically on a K3 surface;\\
$\bullet$ $H$ is a subgroup of $G$;\\
$\bullet$ the K3 surfaces admitting $H$ as group of symplectic
automorphisms admits necessary the group $G$ as group of
symplectic automorphisms.\end{assumeP} In particular we analyzed
this property by a geometrical point of view, constructing
elliptic surfaces and using them to show that certain pairs
$(H,G)$ satisfy property \ref{property P}. Here we reconsider this
property from a point of view more related to the lattices.
\begin{prop}\label{prop: pairs of group (H,G)} A pair of groups $(H,G)$ satisfy property \ref{property P} if
and only if $H\subset G$ and rk$\Omega_H$=rk$\Omega_G$.\end{prop}
\bprf If \ref{property P} holds, then $\Omega_H\subseteq\Omega_G$
(because $H\subset G$ and the K3 surfaces which admits $G$ as
group of symplectic automorphisms admits $H$ as group of
symplectic automorphisms), in particular
rk$\Omega_H\leq$rk$\Omega_G$. If $X$ is a generic algebraic K3
surface such that $H$ acts symplectically on $X$, then
$\rho(X)=1+rk\Omega_H$. If Property \ref{property P} holds, then
$\Omega_{K}\subset NS(X)$ (c.f. Proposition \ref{prop: X has G as
symplectic group iff omegaG in NS(X)}). Since $\Omega_G$ is
negative definite,
$rk\Omega_G\leq \rho(X)-1=rk\Omega_H$, hence rk$\Omega_H=$rk$\Omega_G$.\\
Viceversa if $H\subset G$, then $\Omega_H\subseteq \Omega_G$.
Since rk$\Omega_H$=rk$\Omega_G$, the inclusion
$\Omega_H\hookrightarrow \Omega_G$ has a finite index. As both
$\Omega_H$ and $\Omega_G$ are primitively embedded in the
N\'eron--Severi group of the surface the index of the inclusion
$\Omega_H\hookrightarrow \Omega_G$ is 1, which implies Property
\ref{property P}.\eprf We recall that for each group $G$ of
symplectic automorphisms on a K3 surface rk$\Omega_G$=rk$M_G$ and
rk$M_G$ is computed by \cite{Xiao}.
\begin{rem}{\rm  There are no pairs $(H,G)$ with the property \ref{property P} and
rk$\Omega_G< 16$. Indeed from Xiao's list we deduce that
rk$\Omega_G\in\{8,12,14,15\}$. If $G$ is such that rk$\Omega_G=8$,
then $G=\Z/2\Z$. In the other case there is more then one group,
but there are no inclusions among them. For example the groups $G$
such that rk$\Omega_G$=12 are $\Z/3\Z$ and $(\Z/2\Z)^2$, but of
course there are no inclusions between them, and hence it is not
possible to have a pair $(H,G)$ with the property \ref{property P}
and rk$\Omega_G$=12. The other cases are similar.}\erem\end{rem}
\begin{rem}{\rm Let us assume rk$\Omega_G=16$. The only pairs of
groups $(H,G)$ satisfying the property \ref{property P} are
$(\Z/5\Z,\Dh_5)$, $(\Z/6\Z,\Dh_6)$, $(\Z/2\Z\times
\Z/4\Z,\Z/2\Z\times \Dh_4)$, $((\Z/3\Z)^2,\Ag_{3,3})$ (i.e. the
pairs studied in Section \ref{section: dihedral groups of
automorphisms on elliptic fibration with torsion}). This follows
by the list in \cite{Xiao}: the only groups $G$ such that
rk$\Omega_G=16$ are $\Z/5\Z$, $\Dh_5$, $\Z/6\Z$, $\Dh_6$,
$\Z/2\Z\times \Z/4\Z$, $\Z/2\Z\times \Dh_4$, $(\Z/3\Z)^2$,
$\Ag_{3,3}$, $\Ag_4$. So one has only to check that all the
possible inclusions among these groups are those listed above. In
particular the group $\Ag_4$ is the only one which does not appear
in a pairs $(H,G)$.}\erem
\end{rem}
\begin{rem}{\rm In Xiao's list there are 31 groups $G$ acting symplectically on a K3
surface such that rk$\Omega_G$=19. By \cite{mukai} there are 11
maximal groups acting symplectically on a K3 surface and for each
of them the rank of the associated lattice $\Omega$ is 19. Hence
there exist pairs $(H,G)$ such that $(H,G)$ satisfy the property
\ref{property P} and rk$\Omega_H$=19. It would be interesting to
find explicit examples to show this fact from a geometrical point
of view.}\erem\end{rem}

\end{document}